\documentclass[11pt,oneside,reqno]{amsart}
\usepackage{amssymb}
\usepackage{amsmath}
\usepackage{amsfonts}
\usepackage{mathtools}
\usepackage{color}
\usepackage{tikz-cd}
\usetikzlibrary{decorations.markings}

\usepackage{xcolor}
\definecolor{deepgreen}{cmyk}{0.99998,0,1,0}

\usepackage{relsize}

\usepackage{hyperref}
\hypersetup{hypertex=true,
	colorlinks=true,
	linkcolor=blue,
	anchorcolor=blue,
	citecolor=blue}
\usepackage{comment}
\usepackage{mathrsfs}

\allowdisplaybreaks[4]

\usepackage{tabularx}
\usepackage{dsfont}
\usepackage{enumitem}

\usepackage{etoolbox}

\usepackage{bm}
\usepackage{bbm}

\usepackage[new]{old-arrows}

\theoremstyle{definition}
\newtheorem{defi}{$\mathbf{Definition}$}[section]
\newtheorem*{pro}{$\mathbf{Proof}$}

\theoremstyle{plain}
\newtheorem{theo}[defi]{$\mathbf{Theorem}$}
\newtheorem{lemma}[defi]{$\mathbf{Lemma}$}

\newtheorem{prop}[defi]{$\mathbf{Proposition}$}

\DeclareMathOperator{\op}{Op}

\DeclareMathOperator{\Sp}{Sp}

\DeclareMathOperator{\Tr}{Tr}

\newcommand{\1}{\mathds{1}}

\newcommand{\pa}{\partial}

\newcommand{\lv}{\left\vert}
\newcommand{\rv}{\right\vert}
\newcommand{\lV}{\left\Vert}
\newcommand{\rV}{\right\Vert}

\newcommand{\blk}{\big(}
\newcommand{\brk}{\big)}

\newcommand{\Blk}{\Big(}
\newcommand{\Brk}{\Big)}

\newcommand{\bv}{\big\vert}
\newcommand{\Bv}{\Big\vert}
\newcommand{\bbv}{\bigg\vert}

\newcommand{\bV}{\big\Vert}
\newcommand{\BV}{\Big\Vert}
\newcommand{\bbV}{\bigg\Vert}

\newcommand{\C}{\mathbb{C}}

\newcommand{\N}{\mathbb{N}}

\newcommand{\R}{\mathbb{R}}

\newcommand{\bT}{\mathbb{T}}

\newcommand{\Z}{\mathbb{Z}}

\newcommand{\cH}{\mathcal{H}}

\newcommand{\rS}{\mathscr{S}}

\makeatletter
\newcommand{\itemtaglabel}[2]{\protected@edef\@currentlabel{[#1]}%
	\label{#2}%
}
\makeatother

\definecolor{pink}{RGB}{249,164,186}
\definecolor{grassgreen}{RGB}{128,255,0}

\numberwithin{equation}{section}

\title{Discrete Mixed Quantization for Cat Maps}

\author{Elena Kim$^1$}
\address{$^{1}$Harvard University, Center of Mathematical Sciences and Applications}
\email{elenakim@fas.harvard.edu}

\author{Qiaochu Ma$^{2}$}
\address{$^2$Department of Mathematics, Texas A\&M University}
\email{qiaochu@tamu.edu}

\date{}

\makeatletter
\newcommand*{\rom}[1]{\expandafter\@slowromancap\romannumeral #1@}
\makeatother

\begin{document}

	\begin{abstract}
        We develop a mixed quantization technique for a vector bundle over a fixed regular graph whose fibers are the spaces of quantum states of the torus. Under natural semiclassical and dynamical assumptions, we study the eigenvalues and eigensections of the vector bundle Laplacian, given by the standard graph Laplacian twisted by quantum cat maps. In the semiclassical limit, the eigenvalue density converges to the Kesten--McKay distribution. The eigensections satisfy a quantum ergodicity result. 
	\end{abstract}
	\clearpage\maketitle

	\section{Introduction}

	The quantum cat map is a classical model in quantum chaos. We construct a new model that couples a finite regular graph with quantum cat maps. We study the eigenvalue distribution and eigensection equidistribution of the associated operators.

	\subsection{Main results}

	We first describe the setting and state our main results.
	
	\subsubsection{Quantum cat maps}
	
	Let $\mathrm{Sp}_{2n}(\mathbb{Z})$ denote the group of integer-valued symplectic matrices, i.e., the integer valued matrices that preserve the standard symplectic form. Such matrices are linear automorphisms of  the torus $\mathbb{T}^{2n} \coloneqq \R^{2n} / \Z^{2n}$. For each $N\in2\mathbb{N}$, quantization of the torus gives a Hilbert space $\cH_N$ of quantum states with $\dim_{\mathbb{C}}\cH_N=N^n$. Furthermore, the quantization of the torus assigns to each $f\in C^\infty(\mathbb{T}^{2n})$ an operator $\op_N(f)\in\mathrm{End}(\cH_N)$, known as a \emph{quantum observable}. Moreover, each $E\in\mathrm{Sp}_{2n}(\mathbb{Z})$ is associated with a unitary operator $M_{E,N}\in\mathrm{U}(\cH_N)$, called a \emph{quantum cat map} defined to satisfy the exact Egorov identity
	\begin{equation}
		M_{E,N}^{-1}\op_N(f)M_{E,N}=\op_N(f\circ E),\quad \text{for all } f \in C^\infty(\bT^{2n}).
	\end{equation}
	The operator $M_{E,N}$ is determined only up to a unit complex scalar. Thus, the symplectic group acts projectively on $\cH_N$. Note that we can write both $M_{E,N}$ and $\op_N(f)$ as $N^n \times N^n$ matrices. We study quantum cat maps as $N \rightarrow \infty$, which is equivalent to taking the semiclassical limit.

	\subsubsection{A sequence of graph vector bundles}\label{subsection:graph-vector-bundles}

	Let $X$ be a finite connected $d$-regular graph ($d\geqslant2$). To each oriented edge $(v,v')$ of $X$, assign a matrix $\chi_{v,v'}\in\mathrm{Sp}_{2n}(\mathbb{Z})$ such that $\chi_{v',v}=\chi_{v,v'}^{-1}$.

	For each $N\in2\mathbb{N}$, we note that the unitary quantum cat map $M_{\chi_{v,v'},N}$ associated with $\chi_{v,v'}$ satisfies $M_{\chi_{v,v'},N}=M_{\chi_{v',v},N}^{-1}$ for the appropriate choice of phase constant. Let $C(X,\cH_N)$ denote the space of $\cH_N$-valued functions on the vertices of $X$. The resulting \emph{vector bundle Laplacian} is given by
	\begin{equation}\label{1.8}
		\Delta^{N}\colon C\big(X,\mathcal{H}_N\big)\to C\big(X,\mathcal{H}_N\big),\qquad(\Delta^{N} u)(v)=\sum_{v'\sim v}M_{\chi_{v,v'},N} u(v'),
	\end{equation}
	where $v\sim v'$ means that $v$ and $v'$ are adjacent. The construction of $\Delta^N$ using inverse operators on reverse edges makes $\Delta^{N}$ selfadjoint. Thus, the eigenvalues of $\Delta^{N}$ are real. Throughout this paper, we allow $X$ to have self-loops and multi-edges, with sums over neighbors understood as sums over oriented edges so that distinct edges are counted separately.
	
	We list the eigenvalues of $\Delta^{N}$, counted with multiplicity, and choose an associated orthonormal basis of eigensections
	\begin{equation}\label{a1}
		\Delta^{{N}}u_{N,i}=\lambda_{N,i}u_{N,i},\qquad\lV u_{N,i}\rV^2_{L^2(X,\mathcal{H}_N)}=1,
	\end{equation}
	for $1\leqslant i\leqslant \lv X\rv N^n$, where $\lv X\rv$ denotes the number of vertices of $X$.

	\subsubsection{Two assumptions}\label{subsection:two-assumptions}

	Fix a base vertex of $X$. For a loop $\gamma=(v_k,\ldots,v_0)$ based at this vertex, with $v_k=v_0$, let
	\begin{equation}\label{eq:chi-gamma}
		\chi(\gamma) \coloneqq \chi_{v_k,v_{k-1}}\cdots\chi_{v_1,v_0}.
	\end{equation}
	This matrix depends only on the homotopy class of $\gamma$ and defines a representation of the fundamental group of $X$. We write $\chi(\gamma)$ in terms of square block matrices
	\begin{equation}
		\chi(\gamma)=\begin{pmatrix}
			A_\gamma&B_\gamma\\
			C_\gamma&D_\gamma
		\end{pmatrix}.
	\end{equation}
	
	We say that this representation is \emph{blockwise invertible}, denoted by \textbf{[BIV]}, if for every nonidentity element $\gamma$ of the fundamental group, we have $\det B_\gamma\neq0$. We say that this representation is \emph{ergodic}, denoted by \textbf{[ERG]}, if every function in $L^2(\mathbb{T}^{2n})$ invariant under all $\chi(\gamma)$ is constant almost everywhere. That is, if $f(\chi(\gamma)z)=f(z)$ for every $\gamma\in\Gamma$ and almost every $z\in\mathbb{T}^{2n}$, then $f$ is constant almost everywhere. In \S\,\ref{examples}, we construct some examples satisfying both \textbf{[BIV]} and \textbf{[ERG]}.

	These two assumptions play complementary roles. \textbf{[BIV]} allows us to carry out the semiclassical analysis, while \textbf{[ERG]} is the dynamical assumption corresponding to the ergodicity of the geodesic flow in the usual quantum ergodicity setting.
	
	\subsubsection{Eigenvalue density}

We define the Kesten--McKay distribution \cite{MR109367,MR629617} by
	\begin{equation}\label{1.21n}
		d\mu_{\mathrm{KM}}(\lambda) \coloneqq \frac{d\sqrt{4(d-1)-\lambda^2}}{2\pi\big(d^2-\lambda^2\big)}\mathbbm{1}_{[-2\sqrt{d-1},2\sqrt{d-1}]}(\lambda)d\lambda.
	\end{equation}
	Here, the $d$'s in $d\mu_{\mathrm{KM}}(\lambda)$ and $d\lambda$ denote measure notation, and should not be confused with the degree $d$ of the graph. 
	
	This is a classical spectral limit for large $d$-regular graphs that become locally tree like, where short cycles contribute negligibly. Although our $X$ remains fixed, \textbf{[BIV]} has an analogous effect as $N\to\infty$. We now state the resulting eigenvalue density law.
	
	\begin{theo}\label{t1.2}
		Assume \textbf{[BIV]}. Then for every interval $I\subseteq \mathbb{R}$,
		\begin{equation}\label{1.20n}
			\lim_{\substack{N\in2\mathbb{N},\\ N\to\infty}}\lv\frac{1}{\lv X\rv\cdot N^{n}}\lv\{i\mid\lambda_{N,i}\in I\}\rv-\int_{I}d\mu_{\mathrm{KM}}(\lambda)\rv=0.
		\end{equation}
	\end{theo}
	We prove Theorem~\ref{t1.2} in \S\,\ref{s5.2}.

	\subsubsection{Eigensection equidistribution}
	
	For $f\in C^\infty(X\times\mathbb{T}^{2n})$, quantizing $f(v,\cdot)$ at each vertex defines an operator on $C(X,\mathcal{H}_N)$ by setting
	\begin{equation}
		\big(\op_N(f)u\big)(v)=\op_N\big(f(v,\cdot)\big)u(v).
	\end{equation}
	
	By Born's rule, the matrix element $\langle\op_N(f)u_{N,i},u_{N,i}\rangle$ is the expectation of the observable $f$ in the eigensection $u_{N,i}$. Such expectations detect how the eigensection is distributed over the graph and torus fibers. Hence, when they approach the uniform average of $f$ for every smooth observable, the eigensections tend to be equidistributed. We now state a quantum ergodicity result showing that this is the case for most eigensections as $N\to\infty$.
	
	\begin{theo}\label{C9'}
		Assume \textbf{[BIV]} and \textbf{[ERG]}. Then for every $f\in C^\infty(X\times \mathbb{T}^{2n})$,
		\begin{equation}\label{9.}
			\begin{split}
				\lim_{\substack{N\in2\mathbb{N},\\ N\to\infty}}\frac{1}{\lv X\rv\cdot N^{n}}\sum_{i}\Bv\big\langle \op_{N}(f)u_{N,i},u_{N,i}\big\rangle_{L^2(X,\mathcal{H}_N)}&\\
				-\frac{1}{\lv X\rv}\sum_{v\in X}\int_{\mathbb{T}^{2n}}f(v,z)dz&\Bv^2=0.
			\end{split}
		\end{equation}
	\end{theo}
	By a standard diagonal argument, we can obtain a density-one sequence of eigensections that tend to be equidistributed as $N\to\infty$. We refer to	\S\,\ref{s9.2x} for the full statement and the proof of Theorem \ref{C9'}.

	\subsection{Techniques and related results}\label{subsection:related-results}

	We now review the results related to our main theorems.
	
	\subsubsection{Mixed quantization}

	The main difficulty in proving our results is simultaneously analyzing the graph and torus variables. To overcome this difficulty, we develop a discrete mixed quantization framework that combines semiclassical analysis on the base graph with geometric quantization on the torus fibers. The key difference from earlier mixed quantization results is that the holonomy in earlier results comes from representations of \emph{compact Lie groups}, whereas ours comes from a representation into the \emph{noncompact group} $\mathrm{Sp}_{2n}(\mathbb Z)$. Therefore, in the present setting, the underlying classical dynamics can exhibit hyperbolic behavior, which is an important feature and calls for different estimates. In particular, the need to control this behavior leads us to identify a new \emph{pseudolocality property} of quantum cat maps, which plays a central role in the proofs of our main results (see \S\S\,\ref{s3.5} and \ref{s7.3}).

	Mixed quantization was first developed in heat kernel analysis by Bismut--Ma--Zhang~\cite{MR2838248,MR3615411}, Ma~\cite{MR4665497}, and Puchol~\cite{MR4611826}. It was subsequently developed for semiclassical analysis of vector bundles over manifolds by Ma--Ma~\cite{MR4808253} and Ben Ovadia--Ma--Rodriguez Hertz~\cite{ovadia2025mixedquantizationpartialhyperbolicity}. More recently, Ma~\cite{disc} introduced mixed quantization for graph vector bundles and applied it to several asymptotic spectral problems.

	\subsubsection{Discrete Quantum Ergodicity}

In this paper, we fix a $d$-regular graph and study the asymptotic behavior of its vector bundle Laplacians as the torus quantization parameter $N\to\infty$. A natural point of comparison is the large graph limit, in which the number of vertices tends to infinity along a sequence of $d$-regular graphs.
	
	The first large graph quantum ergodicity result is due to Anantharaman--Le Masson  ~\cite{AlM15}, who adapted the microlocal proof of quantum ergodicity on manifolds to study $d$-regular expander graphs. Anantharaman ~\cite{MR3649482} later gave several new proofs of this result. Quantum ergodicity for Schr\"odinger operators on graphs was further established by Anantharaman--Sabri ~\cite{MR3961083} for disordered systems and by McKenzie--Sabri ~\cite{MS23} for periodic graphs. Related results for averaging operators on the two-sphere were obtained by Brooks--Le Masson--Lindenstrauss ~\cite{MR3567266}. For a more detailed literature review, we refer to the survey of Anantharaman--Sabri ~\cite{AS19}.

\subsubsection{Quantum cat maps}

Our vector bundle Laplacian $\Delta^N$ generally involves several quantum cat maps through the holonomy of the underlying bundle. This remains true even when $X$ has just one vertex and several self-loops, since $\Delta^N$ is then a sum of the maps assigned to those loops and their inverses. The behavior of individual quantum cat maps has been studied extensively in the literature.

  Quantum cat maps were introduced by Hannay--Berry~\cite{MR602111}. Quantum cat maps for which the underlying symplectic matrix is hyperbolic are a popular discrete model in quantum chaos, since hyperbolic matrices have expanding and contracting directions, much like the Anosov geodesic flow on a negatively curved manifold. Quantum cat maps thus provide a concrete model for studying how classical hyperbolicity affects quantum eigenstates.

	In a result similar to Theorem~\ref{C9'}, Bouzouina--De Bièvre~\cite{MR1387942} established quantum ergodicity for quantum cat maps, showing that a density-one family of eigenfunctions equidistributes as $N\to\infty$ for the full sequence of $N \in \N$. Kurlberg--Rudnick~\cite{KR01-1} proved a different version for two-dimensional cat maps, in which every eigenfunction equidistributes along a density-one sequence of levels $N$. Kurlberg--Ostafe--Rudnick--Shparlinski~\cite{KORS24} extended the latter result to higher-dimensional cat maps.
	
	Notably, the stronger quantum unique ergodicity property that every sequence of eigenfunctions equidistributes can fail. Counterexamples were shown in Faure--Nonnenmacher--De Bièvre~\cite{MR2000926}, Kelmer~\cite{Kel07}, and Kim ~\cite{Ki24}. Work towards understanding this gap between quantum ergodicity and quantum unique ergodicity includes work on entropy by Brooks ~\cite{Bo10}, Faure--Nonnenmacher ~\cite{FN04}, and Rivière ~\cite{Ri11}, and work on the support of semiclassical measures by Schwartz ~\cite{Sc21}, Dyatlov--Jézéquel ~\cite{DJ23}, and Kim ~\cite{Ki24}.

\subsection{Organization of the paper}
	
The paper is divided into two parts, \S\S\,\ref{s1}-\ref{s5} and  \S\S\,\ref{S6}-\ref{s6n}, devoted to Theorems~\ref{t1.2} and~\ref{C9'}, respectively. Readers interested in Theorem~\ref{t1.2} can review the first part independently of the second. 

In the first part, we construct graph vector bundles from representations of the fundamental group in \S\,\ref{s1}, recall quantum translations, metaplectic transformations, and quantum cat maps in \S\,\ref{s3}, and construct a sequence of graph vector bundles in \S\,\ref{s4.1xx}. We prove Theorem~\ref{t1.2}, the Kesten--McKay law, in \S\,\ref{s5}. 

In the second part, we introduce kernel operators and analyze their commutator dynamics in \S\,\ref{S6}, define the quantization of observables on the torus in \S\,\ref{S7}, and combine the graph and torus constructions to obtain discrete mixed quantization in \S\,\ref{S8}. We prove Theorem~\ref{C9'}, eigensection equidistribution, in \S\,\ref{s6n}. 

Finally, we give some examples to which our main results apply in \S\,\ref{examples}.

\subsection{Acknowledgments} 
	
	Elena Kim is supported by NSF DMS-2602417 and the Center for Mathematical Sciences and Applications at Harvard University. Qiaochu Ma is supported by NSF DMS-2247313.

	\section{Graph Vector Bundles}\label{s1}

In this section, we discuss unitary vector bundles over graphs and give some of their basic properties. In \S\,\ref{s2.1}, we construct unitary vector bundles over graphs from representations of the fundamental group. Although this differs from the construction in \S\,\ref{subsection:graph-vector-bundles}, where unitary matrices are assigned to edges, the two constructions are equivalent in general. We shall make this equivalence explicit in \S\,\ref{s4.1x}. In \S\,\ref{s2.2}, we introduce a trace formula for vector bundle Laplacians.
	
	In this section, we mainly refer to~\cite[\S\,2]{disc}.
	
	\subsection{Vector bundles via fundamental group representations}\label{s2.1}

	Let $X$ be a finite connected $d$-regular graph. Its universal covering $\widetilde{X}$ is isomorphic to $\mathbf{T}_d$, the $d$-regular tree. We denote  $\pi_1(X)$, the fundamental group of $X$,  by $\Gamma$. Then $\Gamma$ is a finitely generated free group which acts freely on $\mathbf{T}_d$ by graph automorphisms. We fix a vertex fundamental domain $D\subset\mathbf{T}_d$ for this action. To summarize:
	\begin{equation}\label{n2.1}
		\big(X,\widetilde{X},\pi_1(X),\text{a fundamental domain}\big)\cong\big(\Gamma\backslash \mathbf{T}_d, \mathbf{T}_d,\Gamma,D\big).
	\end{equation}
	We shall use $v$ and $\bm{v}$, possibly with subscripts, to denote vertices of $X$ and $\mathbf{T}_d$, respectively. We abuse notation and use  $X$ and $\mathbf{T}_d$ to denote both the graphs and their respective sets of vertices. 
	
	We recall a path description of the fundamental group. For $v\in X$, let $\Omega(X,v)$ denote the space of loops based at $v$,
	\begin{equation}\label{n2.2}
		\Omega(X,v)=\big\{\text{paths in $X$ starting and ending at $v$}\big\}.
	\end{equation}
	We define an equivalence relation on $\Omega(X,v)$ by identifying
	\begin{equation}
		(v_j,\cdots,v_{i+1},v_i,v_{i-1},\cdots,v_0) \quad \text{with} \quad (v_j,\cdots,v_{i+2},v_{i-1},\cdots,v_0)
	\end{equation}
	whenever $v_{i-1}=v_{i+1}$. Thus, two loops are equivalent if one can be obtained from the other by repeatedly inserting or removing immediate backtracking steps. The fundamental group (based at $v$) is the quotient
	\begin{equation}\label{x2.4}
		\Gamma=\Omega(X,v)\big/\sim.
	\end{equation}
	We shall omit the dependence on the base point $v$ from the notation.
	
	Let
	\begin{equation}\label{unitary}
		\rho\colon \Gamma\to \mathrm{U}_\ell(\mathbb{C})
	\end{equation}
	be an $\ell$-dimensional unitary representation.  
	The representation $\rho$ induces a unitary vector bundle $F$ over $X$ by
	\begin{equation}\label{n2.4}
		\begin{split}
			F &\coloneqq \Gamma\backslash\big(\mathbf{T}_d\times\mathbb{C}^\ell\big)\\
			&= \big\{(\bm{v},s)\in\mathbf{T}_d\times\mathbb{C}^\ell\big\}/(\bm{v},s)\sim(\gamma \bm{v},\rho(\gamma)s)\ \text{for every }\gamma\in\Gamma.
		\end{split}
	\end{equation}
	Its fiber dimension $\dim_{\mathbb{C}}F_{v}$ at $v\in X$ and its total dimension $\dim_{\mathbb{C}}F$ are defined by 
	\begin{equation}\label{x2.5}
		\dim_{\mathbb{C}}F_{v} \coloneqq \dim_{\mathbb{C}}\mathbb{C}^\ell=\ell,\qquad\dim_{\mathbb{C}}F \coloneqq \sum_{v\in X}\dim_{\mathbb{C}}F_{v}=\lv X\rv\ell.
	\end{equation}
	The space $C(X,F)$ of sections of $F$ is naturally identified with the space of $\Gamma$-equivariant $\mathbb{C}^\ell$-valued functions on $\mathbf{T}_d$, that is,
	\begin{equation}\label{n2.5}
		C(X,F)=\big\{u\colon \mathbf{T}_d\to \mathbb{C}^\ell\mid u(\bm{v})=\rho(\gamma)u(\gamma^{-1}\bm{v})\text{ for every } \bm{v}\in\mathbf{T}_{d}, \gamma\in\Gamma\big\}.
	\end{equation}

	Since $\rho$ is unitary, the standard Hermitian norm $\lV\cdot\rV_{\mathbb{C}^\ell}$ induces an $L^2$-norm $\lV\cdot\rV_{L^2(X,F)}$ on $C(X,F)$,
	\begin{equation}\label{eq:hermitian -induces}
		\lV u\rV_{L^2(X,F)}^2 \coloneqq \sum_{\bm{v}\in D}\lV u(\bm{v})\rV_{\mathbb{C}^\ell}^2.
	\end{equation}

	We define the Laplacian $\Delta^F$, acting on $C(X,F)$, by
	\begin{equation}
		(\Delta^Fu)(\bm{v})=\sum_{\bm{v}'\sim \bm{v}}u(\bm{v}').
	\end{equation}
	Since $\Delta^F$ is selfadjoint with respect to the inner product on $L^2(X,F)$, we can write an orthonormal basis of eigensections and their respective eigenvalues
	$\{\lambda_{F,i}, u_{F,i}\}_{i=1}^{\dim_{\mathbb{C}}F}$, 
	\begin{equation}\label{n2.8}
		\Delta^{F}u_{F,i}=\lambda_{F,i}u_{F,i},\qquad\lV u_{F,i}\rV^2_{L^2(X,F)}=1.
	\end{equation}

	\subsection{A trace formula}\label{s2.2}

	We give an elementary trace formula that will be used in the analysis of eigenvalues.
	
	\begin{lemma}
		Let $F$ be a vector bundle over $X$ of the form~\eqref{n2.4}. Then for $k\in\mathbb{N}$,
		\begin{equation}\label{n2.9}
			\mathrm{Tr}^{F}\big[\big(\Delta^{F}\big)^k\big]=\sum_{i=1}^{\dim_\mathbb{C}F}\lambda_{F,i}^k=\sum_{{v}\in X,\gamma\in \Omega(X,{v}),\lv\gamma\rv=k}\mathrm{Tr}^{F_{{v}}}\big[\rho(\gamma)\big],
		\end{equation}
		where $\Omega(X,{v})$ is the loop space defined in~\eqref{n2.2} and $\lv\gamma\rv$ denotes the length of $\gamma$.
	\end{lemma}
	
\begin{pro}
Expanding $(\Delta^F)^k$ gives a sum over paths of length $k$. For a vertex $v\in X$, only walks that return to $v$ contribute to its diagonal block, and each such walk $\gamma$ contributes the holonomy $\rho(\gamma)$.\qed
\end{pro}

	\section{Quantum Cat Maps}\label{s3}

	In this section, we introduce quantum cat maps, which quantize integer symplectic matrices and play the role of quantum evolution operators. We present only the material needed for Theorem~\ref{t1.2}. The additional background required for Theorem~\ref{C9'} is introduced in \S\,\ref{S7}.  In \S\,\ref{s3.1}, we define quantum translations. In \S\,\ref{s3.2}, we introduce metaplectic transformations. In \S\,\ref{s3.3}, we define the space of quantum states. In \S\,\ref{s3.4}, we define quantum cat maps as metaplectic transformations acting on the space of quantum states. In \S\,\ref{s3.5}, we establish a decay estimate for the normalized traces of quantum cat maps.
	
	In this section, we mainly refer to ~\cite[\S\,11]{MR2952218},~\cite[\S\,2]{DJ23}, and ~\cite[\S\,2]{Ki24}.
	
	\subsection{Quantum translations}\label{s3.1}

	Let $\rS(\R^{n})$ be the Schwartz space of rapidly decreasing functions, with its dual space $\mathscr{S}'(\R^{n})$ of tempered distributions. 
	
	Take a semiclassical parameter $0<h\leqslant 1$. For $\omega=(y,\eta) \in\R^{2n}$, we define the associated \emph{quantum translation} $U_{\omega}$. We give a formal definition in ~\eqref{3.3}, but for now, we use only the following explicit formula: 
	\begin{equation}\label{3.1x}
		\begin{split}
			&U_{\omega}\colon \mathscr{S}(\R^{n})\to \mathscr{S}(\R^{n}),\\
			&(U_{\omega} f)(x)=e^{\frac{\sqrt{-1}}{h}\langle y,x\rangle+\frac{\sqrt{-1}}{2h} \langle y,\eta\rangle} f(x+\eta).
		\end{split}
	\end{equation}

	Using~\eqref{3.1x}, we can verify directly that the quantum translations satisfy the following composition formula
	\begin{equation}\label{3.3x}
		U_\omega U_{\omega'}=e^{\frac{\sqrt{-1}}{2h}\sigma(\omega,\omega')} U_{\omega+ \omega'},
	\end{equation}
	where 
	$\sigma$ is the \emph{standard symplectic form} on $\R^{2n}$ given by
	\begin{equation}\label{3.2x}
		\sigma(\omega,\omega') \coloneqq \langle\eta,  y'\rangle-\langle y, \eta'\rangle, \quad\text{for all }\omega=(y,\eta), \text{ } \omega'=(y',\eta') \in \R^{2n}.
	\end{equation}

	\subsection{Metaplectic transformations}\label{s3.2}

	Let $\Sp_{2n}(\R)$ denote the group of real symplectic $2n\times2n$ matrices preserving the symplectic form $\sigma$ defined in~\eqref{3.2x}, that is, for any $E\in\Sp_{2n}(\R)$ and $\omega,\omega'\in\R^{2n}$, we have 
	\begin{equation}\label{eq:symplectic-matrix}
		\sigma(E\omega,E\omega')=\sigma(\omega,\omega').
	\end{equation}
	Write
	\begin{equation}\label{3.5x}
		E=\begin{pmatrix}
			A&B\\
			C&D
		\end{pmatrix}\in\Sp_{2n}(\mathbb{R}), \quad A, B, C, D \in M_n(\R).
	\end{equation}
	We refer to $A$, $B$, $C$, and $D$ as the $A$-, $B$-, $C$-, and $D$-blocks of $E$, respectively. From~\eqref{eq:symplectic-matrix}, we have that
	\begin{equation}\label{3.5xx}
		A^{-1}B=(A^{-1}B)^T,\qquad
		A^TD-C^TB=I.
	\end{equation}
	
	For each $E \in \Sp_{2n}(\R)$, let $\mathcal{M}_E$ be the set of all unitary operators $M_E: L^2(\R^n) \rightarrow L^2(\R^n)$ that satisfy an exact Egorov's theorem~\eqref{MA}. We delay the full discussion of the exact Egorov's theorem~\eqref{MA} until \S~\ref{S7}, but note here a particular case involving the quantum translations~\eqref{3.1x}: 
	\begin{equation}\label{3.8xx}
		M_E^{-1}U_\omega M_E=U_{E^T\omega}, \quad \text{for all } \omega\in\mathbb{R}^{2n}, \text{ } E\in \Sp_{2n}(\R), \text{ } M_E \in \mathcal{M}_E.
	\end{equation}
	We call such $M_E$ \emph{metaplectic transformations} associated with $E$. In particular, such operators are unique up to multiplication by a unit complex scalar. Fixing a particular choice of $M_E$ for each matrix $E$, we have the following commutation formula for any $E,E'\in\mathrm{Sp}_{2n}(\mathbb{R})$,
	\begin{equation}\label{3.6x}
		M_EM_{E'}=c(E,E')M_{EE'},
	\end{equation}
	where $c(E,E')$ is a unit complex factor. The results in this paper are invariant under the choice of $M_E$.

	By ~\cite[Theorem 11.10]{MR2952218}, if the $A$-block in~\eqref{3.5x} satisfies $\det A \neq 0$, then up to a unit factor, we have the following explicit formula for $M_E$, 
	\begin{equation}\label{3.8x}
		\begin{split}
			M_E&\colon \mathscr{S}(\mathbb{R}^n)\to \mathscr{S}(\mathbb{R}^n),\\
			(M_Ef)(x)&=\frac{\lv\det A\rv^{-\frac{1}{2}}}{(2 \pi h)^n} \int_{\R^n} \int_{\R^n} e^{\frac{\sqrt{-1}}{h} (\varphi(x,{\eta}) -\langle {y},{\eta}\rangle )}f({y}) d{y}d{\eta},\\
			\varphi(x,{\eta})&\coloneqq \frac{1}{2}\langle CA^{-1}x,x\rangle+\langle A^{-1}x, {\eta} \rangle - \frac{1}{2} \langle A^{-1} B{\eta},{\eta} \rangle.
		\end{split}
	\end{equation}
	A general $E\in\mathrm{Sp}_{2n}(\mathbb{R})$ can be factored as $E=E_1E_2$, where the $A$-blocks of both $E_1$ and $E_2$ are invertible. Hence the general metaplectic transformation can be obtained, up to a unit complex scalar, by composing two transformations of the form~\eqref{3.8x}.

	\subsection{The space of quantum states}\label{s3.3}

	For the rest of the paper, we assume that the semiclassical parameter in~\eqref{3.1x} and~\eqref{3.8x} is given by
	\begin{equation}\label{x3.9}
		h= (2\pi N)^{-1},\quad N\in 2\mathbb{N}.
	\end{equation}
	We define \emph{the space of quantum states} by
	\begin{equation}\label{3.11f}
		\cH_N\coloneqq \big\{f \in \rS'(\R^n)\mid U_\omega f=e^{N \pi \sqrt{-1}\langle y,\eta \rangle}f  \text{ for all } \omega=(y,\eta) \in \Z^{2n}\big\}.
	\end{equation}
	By \cite[Proposition 2.1]{MR1387942}, our choice of $h$ ensures that the spaces of quantum states are nontrivial. 
	We note that ~\cite[\S\,2.2.1]{DJ23} gives a more general definition for quantum states with an additional parameter $\theta \in \bT^{2n}$. In our paper, the space $\cH_N$ corresponds to the case where $\theta=0$.

	Define
	\begin{equation}
		\Z^n_N \coloneqq (\mathbb Z \bmod{N})^n\cong\{0,\cdots,N-1\}^n.
	\end{equation}
	We have the following explicit basis of $\mathcal{H}_N$.
	\begin{lemma}[{\cite[Lemma 2.5]{DJ23}}]
		The space $\cH_N$ has dimension $N^n$ and admits a basis $\{e_j\}_{j\in\mathbb{Z}_N^n}$ given by
		\begin{equation}\label{3.12x}
			e_j(x) \coloneqq N^{-\frac{n}{2}}\sum_{k \in \Z^n} \delta\big(x-\tfrac{Nk+j}{N}\big).
		\end{equation}
	\end{lemma}
	Let $\langle \cdot, \cdot\rangle_{\cH_N}$ be the inner product on $\cH_N$ such that $\{e_j\}$ is orthonormal.

	We shall see that certain quantum translations and all metaplectic transformations act on $\mathcal{H}_N$. The action of quantum translations follows directly from their definition and is described below, while the action of metaplectic transformations will be considered in the next subsection.

	\begin{lemma}
		For $\omega=(y,\eta)\in \mathbb{Z}^{2n}$, the quantum translation $U_{\omega/N}$ acts on $\mathcal{H}_N$. Indeed, for $j\in \mathbb{Z}_N^n$,
		\begin{equation}\label{3.19x}
			U_{\omega/N}e_j=e^{\frac{2 \pi\sqrt{-1} }{N}\langle y, j \rangle -\frac{\pi\sqrt{-1} }{N} \langle y,\eta\rangle} e_{j-\eta}.
		\end{equation}
		
	\end{lemma}

	\begin{pro}
		By~\eqref{3.1x}, we have
		\begin{equation}
			\begin{split}
				&\big(U_{\omega/N}e_j\big)(x)\\
				&=e^{2\pi\sqrt{-1}\langle y, x \rangle+\frac{\pi\sqrt{-1}}{N}\langle y,\eta\rangle }N^{-\frac{n}{2}} \sum_{k \in \Z^n} \delta\big(x+\tfrac{\eta}{N}-\tfrac{Nk+j}{N}\big)\\
				&=e^{\frac{2\pi\sqrt{-1}}{N}\langle  y,j-\eta\rangle+\frac{\pi\sqrt{-1}}{N}\langle y,\eta\rangle }N^{-\frac{n}{2}} \sum_{k \in \Z^n} \delta\big(x-\tfrac{Nk+j-\eta}{N}\big),
			\end{split}
		\end{equation}
		which gives~\eqref{3.19x}.\qed
	\end{pro}

	\subsection{Quantum cat maps}\label{s3.4}

	We assumed that $N$ is even in~\eqref{x3.9} so that following ~\cite[\S\,2.2.4]{DJ23}, for any $E\in\mathrm{Sp}_{2n}(\mathbb{Z})$, the metaplectic transform $M_E$  acts on $\mathcal{H}_N$. We denote this by 
	\begin{equation}
		M_{E,N}=M_E\big|_{\mathcal{H}_N},
	\end{equation}
	and call it the \emph{quantum cat map} associated with $E$ of level $N$.

	We now use~\eqref{3.8x} to compute an explicit formula for $M_{E,N}$ with respect to the basis~\eqref{3.12x}. On a first reading, readers may accept the formula and return to its proof later.
	
	\begin{prop}
		Let $E\in \Sp_{2n}(\Z)$ as in~\eqref{3.5x}. Assume that $\det A\neq0,\det B\neq0$. Then for $j,j'\in \mathbb{Z}_N^n$,
		\begin{equation}\label{n2.8'x}
			\begin{split}
				&\langle M_{E,N} e_j, e_{j'} \rangle_{\cH_N}\\
				&= \frac{e^{-\sqrt{-1}\frac{\pi}{4} \mathrm{sgn}(A^{-1}B)}}{|N^n\det B|^{\frac{1}{2}}} e^{\frac{\pi\sqrt{-1}}{N} (\langle DB^{-1}{j'},{j'}\rangle - 2\langle B^{-1}{j'}, j \rangle + \langle B^{-1} A j, j\rangle )}\\
				&\quad\cdot\sum_{k \in \Z^n/(B^T\Z^n)}  e^{\pi\sqrt{-1} (N \langle B^{-1} A k,k \rangle + 2\langle B^{-1} A k, j\rangle - 2 \langle B^{-1}{j'}, k \rangle)}.
			\end{split}
		\end{equation}
	\end{prop}
	
	\begin{pro}
		
		Let $g(\eta) \in C_c^\infty(\R^n)$ such that $g({\eta}) =1$ for $\lv{\eta}\rv\leqslant 1$ and $g({\eta}) =0$ for $\lv{\eta}\rv > 2$. We can extend~\eqref{3.8x} to $f= e_j$, where the ${y}$-integral is understood distributionally against the Dirac comb. Then by the Poisson summation formula
		\begin{equation}\label{3.28'x}
			\begin{split}
				&\big(M_{E,N} e_j\big)(x)\\
				&=\lim_{\varepsilon \rightarrow 0} \frac{N^{\frac{n}{2}}}{\lv\det A\rv^{\frac{1}{2}}}\int_{\R^{2n}}g(\varepsilon {\eta})e^{\frac{\sqrt{-1}}{h} (\varphi(x, {\eta}) -\langle {y}, {\eta} \rangle )}\sum_{k \in \Z^n} \delta({y}-k-\tfrac{j}{N}) d{y}d{\eta}\\
				&=\lim_{\varepsilon \rightarrow 0}  \frac{N^{\frac{n}{2}}}{|\det A|^{\frac{1}{2}}}\int_{\R^n} g(\varepsilon {\eta}) e^{\frac{\sqrt{-1}}{h}( \varphi(x,{\eta}) -\langle \frac{j}{N},{\eta}\rangle)}\sum_{k \in \Z^n} e^{-\frac{\sqrt{-1}}{h}  \langle k,{\eta} \rangle}d{\eta}.
			\end{split}
		\end{equation}
		Evaluating the oscillatory Gaussian integral in~\eqref{3.28'x} and using~\eqref{3.5xx}, we obtain
		\begin{equation}\label{3.29x}
			\begin{split}
				(M_{E,N}e_j)(x) &=\frac{N^{\frac{n}{2}}}{|\det A|^{\frac{1}{2}}} e^{\frac{\sqrt{-1}}{2h} \langle CA^{-1}x,x\rangle}\\
				&\cdot\sum_{k \in \Z^n} \lim_{\varepsilon \rightarrow 0}  \int_{\R^n} g(\varepsilon {\eta}) e^{\frac{\sqrt{-1}}{h}(\langle A^{-1}x-k -\frac{j}{N}, {\eta}\rangle-\frac{1}{2}\langle A^{-1}B {\eta},{\eta}\rangle)}d{\eta}\\
				&=\frac{1}{\vert\det B\vert^{\frac{1}{2}}}e^{-\frac{\pi\sqrt{-1}}{4} \mathrm{sgn}(A^{-1}B)+\frac{\sqrt{-1}}{2h} \langle CA^{-1}x,x\rangle}\\
				&\quad\cdot\sum_{k\in\Z^n}e^{\frac{\sqrt{-1}}{2h}{\langle B^{-1}x-B^{-1}A(k+\frac{j}{N}), A^{-1}x-(k+\frac{j}{N})}\rangle}\\
				&=\frac{1}{|\det B|^{\frac{1}{2}}}e^{-\frac{\pi\sqrt{-1}}{4} \mathrm{sgn}(A^{-1}B)+\frac{\sqrt{-1}}{2h} \langle CA^{-1}x,x\rangle}\\
				&\quad\cdot\sum_{k\in\Z^n}e^{\frac{\sqrt{-1}}{2h}\langle B^{-1}x,A^{-1}x\rangle
					-\frac{\sqrt{-1}}{h}\langle B^{-1}x,(k+\frac{j}{N})\rangle
					+\frac{\sqrt{-1}}{2h}\langle B^{-1}A(k+\frac{j}{N}),(k+\frac{j}{N})\rangle}.
			\end{split}
		\end{equation}
		Applying~\eqref{3.5xx} again, we get $DB^{-1}=(A^{T})^{-1}C^T+(A^{T})^{-1}B^{-1}$. Using this to simplify the phase in~\eqref{3.29x}, we get
		\begin{equation}\label{3.29xx}
			\begin{split}
				(M_{E,N}e_j)(x) =&\frac{1}{\lv\det B\rv^{\frac{1}{2}}}e^{-\frac{\pi\sqrt{-1}}{4} \mathrm{sgn}(A^{-1}B)+\frac{\sqrt{-1}}{2h} \langle DB^{-1}x,x\rangle}\\
				&\cdot\quad\sum_{k \in \Z^n}e^{\frac{\sqrt{-1}}{h}( -  \langle B^{-1}x,k+\frac{j}{N} \rangle + \frac{1}{2}\langle B^{-1}A(k+\frac{j}{N}),k+\frac{j}{N}\rangle)}.
			\end{split}
		\end{equation}
		Then using $\det B\neq0$,  ~\eqref{x3.9}, and the Poisson summation formula, we have 
		\begin{equation}\label{3.17xx}
			\begin{split}
				&\sum_{k \in \Z^n}e^{\frac{\sqrt{-1}}{h}( -  \langle B^{-1}x,k+\frac{j}{N} \rangle + \frac{1}{2}\langle B^{-1}A(k+\frac{j}{N}),k+\frac{j}{N}\rangle)}\\
				&=\sum_{\substack{k \in \Z^n/(B^T\Z^n),\\ k'\in \Z^n}}e^{\frac{\sqrt{-1}}{h}( -\langle B^{-1}x,k+B^T{k'} +\frac{j}{N} \rangle + \frac{1}{2}\langle   B^{-1} A (k+B^T{k'}+\frac{j}{N}),  k+B^T{k'} +\frac{j}{N}\rangle)}\\
				&=\sum_{k \in \Z^n/(B^T\Z^n)}  e^{2 \pi N \sqrt{-1}\langle B^{-1}(-x+ A(\frac{k}{2}+\frac{j}{2N})), k +\frac{j}{N}\rangle} \sum_{{k'} \in \Z^n} e^{2 \pi N \sqrt{-1} \langle-x, {k'} \rangle}\\
				&=\frac{1}{N^n} \sum_{k \in \Z^n/(B^T\Z^n)}  e^{2 \pi N\sqrt{-1}\langle B^{-1}(-x+A(\frac{k}{2}+\frac{j}{2N})), k +\frac{j}{N}\rangle} \sum_{{k'} \in \Z^n} \delta\big(x+\tfrac{{k'}}{N}\big).
			\end{split}
		\end{equation}

		By~\eqref{3.29xx} and~\eqref{3.17xx}, we conclude
		\begin{equation}\label{n2.8'}
			\begin{split}
				\langle M_{E,N} e_j, e_{j'} \rangle_{\cH_N}&= \frac{e^{-\sqrt{-1}\frac{\pi}{4} \mathrm{sgn}(A^{-1}B)}}{|N^n\det B|^{\frac{1}{2}}} e^{\frac{\pi\sqrt{-1}}{N} \langle DB^{-1} {j'},{j'} \rangle}\\
				&\quad\cdot\sum_{k \in \Z^n/(B^T\Z^n)}e^{2\pi N\sqrt{-1}\langle B^{-1}(-\frac{{j'}}{N}+ A(\frac{k}{2}+\frac{j}{2N})), k +\frac{j}{N}\rangle},
			\end{split}
		\end{equation}
		which proves~\eqref{n2.8'x}.\qed
	\end{pro}

	\subsection{Trace estimate}\label{s3.5}
	
	We now establish a trace estimate for quantum cat maps, which we view as the first step towards proving pseudolocality in Proposition~\ref{prop:pseudolocality}.

	\begin{prop}\label{p3.2}
		Let $E\in \Sp_{2n}(\Z)$ as in~\eqref{3.5x}. Assume that $\det B\neq0$. Then for $N\in 2\mathbb{N}$ and $j,j'\in \mathbb{Z}_N^n$,
		\begin{equation}\label{3.28x'}
			\lv\langle M_{E,N} e_j, e_{j'} \rangle_{\cH_N}\rv\leqslant\frac{ \lv\det B\rv^{\frac{1}{2}}}{N^{\frac{n}{2}}}.
		\end{equation}
		In particular,
		\begin{equation}\label{3.28x}
			\frac{1}{\dim_\C(\cH_N)}\lv\Tr^{\cH_N}[M_{E,N}]\rv\leqslant \frac{\lv \det B\rv^{\frac{1}{2}}}{N^{\frac{n}{2}}}.
		\end{equation}
	\end{prop}

	\begin{pro}
		If the $A$-block of $E$ is  invertible, then~\eqref{3.28x'} follows directly from~\eqref{n2.8'x}, since
		\begin{equation}\label{3.32x}
			\lv\langle M_{E,N} e_j, e_{j'} \rangle_{\cH_N}\rv\leqslant \frac{\lv\Z^n/(B^T\Z^n)\rv}{\lv N^n\det B\rv^{\frac{1}{2}}}  = \frac{ \lv\det B\rv^{\frac{1}{2}}}{N^{\frac{n}{2}}}.
		\end{equation}
		
		We can also prove~\eqref{3.28x'} without using the explicit formula~\eqref{n2.8'x}. For $y\in\mathbb{Z}_N^n$, by~\eqref{3.8xx}, we get
		\begin{equation}\label{3.33x}
			U_{(y,0)/N}M_{E,N}=M_{E,N} U_{E^T(y,0)/N}=M_{E,N}U_{(A^Ty,B^Ty)/N}.
		\end{equation}
		From~\eqref{3.19x}, we see that $U_{(y,0)/N}$ acts as a diagonal unitary matrix with respect to $\{e_j\}$. Hence
		\begin{equation}\label{3.34x}
			\lv\langle U_{(y,0)/N}M_{E,N}e_j,e_{j'}\rangle\rv=\lv\langle M_{E,N}e_j,e_{j'}\rangle\rv.
		\end{equation}
		Applying~\eqref{3.19x} again, we have
		\begin{equation}\label{3.35'x}
			\lv\langle M_{E,N}U_{(A^Ty,B^Ty)/N}e_j,e_{j'}\rangle\rv=\lv\langle M_{E,N}e_{j-B^Ty},e_{j'}\rangle\rv.
		\end{equation}
		By~\eqref{3.33x},~\eqref{3.34x}, and~\eqref{3.35'x}, for any fixed ${j'}\in \mathbb{Z}_N^n$, the function
		\begin{equation}
			j\in \mathbb{Z}_N^n\mapsto \lv\langle M_{E,N}e_j,e_{j'}\rangle\rv
		\end{equation}
		is constant on cosets of the image of
		\begin{equation}
			B^T\bmod{N}\colon \mathbb{Z}_N^n\to \mathbb{Z}_N^n.
		\end{equation}
		Since $M_{E,N}$ is unitary,
		\begin{equation}\label{3.38'x}
			\begin{split}
				\lv\mathrm{Im}(B^T\bmod{N})\rv\lv\langle M_{E,N}e_j,e_{j'}\rangle\rv^2&=\sum_{j''\in j+\mathrm{Im}(B^T)}\lv\langle M_{E,N}e_{j''},e_{j'}\rangle\rv^2\\
				&\leqslant\sum_{j''\in \mathbb{Z}_N^n}\lv\langle M_{E,N}e_{j''},e_{j'}\rangle\rv^2=1.
			\end{split}
		\end{equation}
		It remains to estimate the cardinality of $\mathrm{Im}(B^T\bmod{N})$. We compute that
		\begin{equation}
			\begin{split}
				&\bv \mathbb{Z}_N^n\big/\mathrm{Im}(B^T\bmod{N})\bv\\
				&=\Bv(\mathbb{Z}^n/(N\mathbb{Z})^n)\big/ \big((B^T\mathbb{Z}^n+(N\mathbb{Z})^n)\big/(N\mathbb{Z})^n\big)\Bv\\
				&=\bv\mathbb{Z}^n\big/(B^T\mathbb{Z}^n+(N\mathbb{Z})^n)\bv\leqslant \bv\mathbb{Z}^n\big/B^T\mathbb{Z}^n\bv\\
				&=\lv\det(B)\rv.
			\end{split}
		\end{equation}
		Hence,
		\begin{equation}\label{3.39x}
			\lv\mathrm{Im}(B^T\bmod{N})\rv\geqslant N^n/\lv\det(B)\rv.
		\end{equation}
		Combining~\eqref{3.38'x} and~\eqref{3.39x}, we get~\eqref{3.28x'}.\qed
	\end{pro}

	\section{A Sequence of Vector Bundles}\label{s4.1xx}

	In this section, we combine graph vector bundles from \S\,\ref{s1} with quantum cat maps from  \S\,\ref{s3} to form the geometric setting of our main results. In \S\,\ref{s4.1x}, we introduce a sequence of vector bundles induced by an integer symplectic representation of the fundamental group of a graph. In \S\,\ref{s4.2}, we introduce two assumptions that play roles analogous to pseudolocality and the ergodicity of the geodesic flow, respectively.
	
	The construction of this sequence of vector bundles from geometric representations is inspired by ~\cite[\S\,4.3]{MR3615411} and ~\cite[\S\,4.2]{disc}.
	
	\subsection{A sequence of vector bundles}\label{s4.1x}
	Let $X=\Gamma\backslash\mathbf{T}_d$ be a $d$-regular graph as in~\eqref{n2.1}.
	Recall that in \S\,\ref{subsection:graph-vector-bundles},  to each oriented edge $(v,v')$ of $X$, we assigned a matrix $\chi_{v,v'}\in\mathrm{Sp}_{2n}(\mathbb{Z})$ such that $\chi_{v',v}=\chi_{v,v'}^{-1}$. Using the definition~\eqref{eq:chi-gamma}, we have a representation  
	\begin{equation}\label{5.7x}
		\chi\colon \Gamma\to\mathrm{Sp}_{2n}(\mathbb{Z}).
	\end{equation}

	By~\eqref{3.6x}, for $N\in2\mathbb{N}$, the representation $\chi$ induces a projective unitary representation of $\Gamma$ on $\mathcal{H}_N$. Since $\Gamma$ is a free group, this projective representation can be lifted to a unitary representation
	\begin{equation}\label{4.2x}
		\rho_N\colon\Gamma\to\mathrm{U}(\mathcal{H}_N),\quad\rho_N(\gamma)=c_N(\gamma)M_{\chi(\gamma),N},
	\end{equation}
	where $c_N(\gamma)$ is a unit complex factor.  Indeed, we fix a set of free generators $(\gamma_i)_i$ of $\Gamma$, choose unitary representatives $\rho_N(\gamma_i)=M_{\chi(\gamma_i),N}$ and extend multiplicatively to all of $\Gamma$.

	We define the associated unitary flat vector bundle $F_N$ by
	\begin{equation}\label{4.5ns}
		\begin{split}
			F_N& \coloneqq \Gamma\backslash\big(\mathbf{T}_d\times\mathcal{H}_N\big)\\
			&=\big\{(\bm{v},s)\in\mathbf{T}_d\times\mathcal{H}_N\big\}/(\bm{v},s)\sim(\gamma\bm{v},\rho_N(\gamma)s)\ \text{for every }\gamma\in\Gamma,
		\end{split}
	\end{equation}
	
	For $F_N$, we record its fiber dimension at $v\in X$, its total dimension, its space of sections, its Laplacian, and its $L^2$-norm:
	\begin{equation}\label{4.5nsss}
		\begin{split}
			\dim_{\mathbb{C}}F_{N,v}
			&\coloneqq \dim_{\mathbb{C}}\mathcal{H}_N,
			\qquad
			\dim_{\mathbb{C}}F_N
			\coloneqq \lv X\rv\dim_{\mathbb{C}}\mathcal{H}_N,\\
			C(X,F_N)
			&\coloneqq \big\{u\colon\mathbf{T}_d\to\mathcal{H}_N\mid
			u(\bm{v})=\rho_N(\gamma)u(\gamma^{-1}\bm{v})\\
			&\qquad\qquad\qquad\qquad\ \ \text{for every }\bm{v}\in\mathbf{T}_d,\ \gamma\in\Gamma\big\},\\
			(\Delta^{F_N}u)(\bm{v})
			&\coloneqq \sum_{\bm{v}'\sim\bm{v}}u(\bm{v}'),\qquad
			\lV u\rV_{L^2(X,F_N)}^2 \coloneqq \sum_{\bm{v}\in D}\lV u(\bm{v})\rV_{\mathcal{H}_N}^2,
		\end{split}
	\end{equation}
	where $D$ is a vertex fundamental domain in $\mathbf{T}_d$. Note that  $\dim_{\mathbb{C}}F_N$ is also the dimension of $C(X,F_N)$.

	Since $\rho_N$ is unitary, the Laplacian $\Delta^{F_N}$ is selfadjoint. Hence, we can write an orthonormal basis of eigensections and their respective eigenvalues
	$(\lambda_{N,i}, u_{N,i})_{i=1}^{\dim_{\mathbb{C}}F_N}$,
	\begin{equation}\label{n2.6}
		\Delta^{F_N}u_{N,i}=\lambda_{N,i}u_{N,i},\qquad\lV u_{N,i}\rV_{L^2(X,F_N)}^2=1.
	\end{equation}
We shall study the asymptotic behavior of these eigenvalues and eigensections as $N\to\infty$. Although here we work with $(C(X,F_N),\Delta^{F_N})$ in \eqref{4.5nsss} rather than $(C(X,\mathcal{H}_N),\Delta^{N})$ in \eqref{1.8}, the following result shows these formulations are equivalent. Thus, it suffices to prove Theorems~\ref{t1.2} and~\ref{C9'} for $\Delta^{F_N}$.

	\begin{prop}\label{lem:unitary-equiv}
	There exists a unitary map $\mathcal U$ making the	following diagram commute,
		\begin{equation}\label{com}
		\begin{tikzcd}[ampersand replacement=\&, column sep=normal,row sep=normal]
			C(X,\mathcal{H}_N)\arrow[r,"\Delta^N"]\arrow[d,"\mathcal{U}"swap]\&C(X,\mathcal{H}_N) \arrow[d,"\mathcal{U}"] \\
			C(X,F_N) \arrow[r,"\Delta^{F_N}"]\& C(X,F_N)
		\end{tikzcd}.
	\end{equation}
	\end{prop}

	\begin{pro}
	Fix a base vertex $v_0\in X$. We identify the vertices of $\mathbf T_d$ with homotopy classes of paths starting at $v_0$. For a path $\bm{v}=(v_i,\cdots,v_0)$, we define
		\begin{equation}\label{4.7}
			M_{\bm{v},N}\coloneq M_{\chi_{v_i,v_{i-1}},N}\cdots M_{\chi_{v_1,v_0},N}.
		\end{equation}
	Similarly, we identify $\Gamma$ with homotopy classes of loops based at $v_0$. For a loop $\gamma=(v_0,v_{-1},\ldots,v_{-j})$ with $v_{-j}=v_0$, the representation in \eqref{4.2x} is given by
		\begin{equation}\label{4.8}
\rho_N(\gamma)=M_{\chi_{v_{0},v_{-1}},N}\cdots M_{\chi_{v_{-j+1},v_{-j}},N}.
		\end{equation}
		
For $u\in C(X,\mathcal H_N)$, we define
		\begin{equation}\label{4.9}
	(\mathcal Uu)(\bm{v})\coloneq M_{\bm{v},N}^{-1}u(v_i).
		\end{equation}
We first verify that $\mathcal Uu\in C(X,F_N)$. The deck action is given by $\gamma\bm{v}=(v_{i},\cdots,v_0=v_{-j},v_{-j+1},\cdots,v_0)$. Hence, by \eqref{4.7}, \eqref{4.8}, and \eqref{4.9},
\begin{equation}
M_{\gamma\bm{v},N}=M_{\bm v,N}\rho_N(\gamma)^{-1},\qquad(\mathcal Uu)(\gamma\bm v)
=\rho_N(\gamma)(\mathcal Uu)(\bm v),
\end{equation}
as required.

Since each $M_{\bm{v},N}$ is unitary and a fundamental domain $D$ contains one lift of each vertex of $X$, the map \eqref{4.9} preserves the $L^2$-norm, and $\mathcal U$ is unitary.

It remains to verify the relation between $\mathcal U$ and the two Laplacians. Every neighbor of $\bm v$ is represented by a path $\bm w=(v_{i+1},v_i,\ldots,v_0)$ with $v_{i+1}\sim v_i$. Therefore,
\begin{equation}
\begin{split}
(\Delta^{F_N}\mathcal Uu)(\bm v)
&=\sum_{\bm w=(v_{i+1},v_i,\cdots,v_0)\sim\bm v}M_{\bm w,N}^{-1}u(v_{i+1})\\
&=M_{\bm v,N}^{-1}\sum_{v_{i+1}\sim v_i}M_{\chi_{v_i,v_{i+1}},N}u(v_{i+1})\\
&=(\mathcal U\Delta^Nu)(\bm v).
\end{split}
\end{equation}
This proves the commutativity of \eqref{com}.\qed
	\end{pro}

	\subsection{Two assumptions}\label{s4.2}
	
	We recall and elaborate on the definitions \textbf{[BIV]} and \textbf{[ERG]} from \S~\ref{subsection:two-assumptions}.
	We say that the representation~\eqref{5.7x} is blockwise invertible, denoted by \textbf{[BIV]}, if for every nonidentity $\gamma\in\Gamma$, we have
	\begin{equation}\label{5.10}
		\det B_{\gamma}\neq0,
	\end{equation}
	where $B_{\gamma}$ denotes the $B$-block of $\chi(\gamma)$, as in~\eqref{3.5x}.
	
	We say that the representation~\eqref{5.7x} is ergodic, denoted by \textbf{[ERG]}, if the $\chi(\Gamma)$-action on $\mathbb{T}^{2n}$ is ergodic, or equivalently, for every $\omega\in\mathbb{Z}^{2n}\setminus\{0\}$, the following orbit set is an infinite set,
	\begin{equation}\label{5.11}
		\big\{\chi(\gamma)^T\omega\mid\gamma\in\Gamma\big\}.
	\end{equation}

	These two assumptions can be compared with those in \cite{AlM15}. In \cite{AlM15}, the underlying graphs converge locally to a regular tree, meaning that short loops become asymptotically negligible. In our setting, the base graph is fixed, so its short loops may be numerous and combinatorially complicated. Nevertheless, as the fiber varies, \textbf{[BIV]}~\eqref{5.10} ensures that the contributions from these short loops become asymptotically negligible, therefore allowing us to carry out the semiclassical analysis. Similarly, \textbf{[ERG]}~\eqref{5.11} plays a role analogous to the expander property in \cite{AlM15}, providing the dynamical condition required for quantum ergodicity.
	
	\section{Kesten--McKay law}\label{s5}

	In \S\,\ref{s5.1}, we define the normalized eigenvalue counting measure of a vector bundle Laplacian. In \S\,\ref{s5.2}, we state and prove our eigenvalue distribution result. For related results on Kesten--McKay law, we refer to ~\cite[Proposition 4.1]{MR1677685},~\cite[Corollary 2]{MR3567266}, and~\cite[Theorem 6.3]{disc}.

	\subsection{Eigenvalue counting measure}\label{s5.1}
	
	In the setting of~\eqref{n2.6}, let $d\mu_{N}(\lambda)$ denote the normalized counting measure of the eigenvalues of $\Delta^{F_N}$, given by
	\begin{equation}\label{4.24}
		d\mu_{N}(\lambda) \coloneqq \frac{1}{\dim_\mathbb{C}F_N}\sum_{i=1}^{\dim_\mathbb{C}F_N}\delta_{\lambda_{N,i}}(\lambda).
	\end{equation}

	We note that the Kesten--McKay distribution defined in~\eqref{1.21n} can be characterized by the following moment identity
	\begin{equation}\label{eq:moment-identity}
		\int_\mathbb{R}\lambda^kd\mu_{\mathrm{KM}}(\lambda)=\bv\{\gamma\in \Omega(\mathbf{T}_d,o)\mid\lv\gamma\rv=k\}\bv
	\end{equation}
	for every $k\in\mathbb{N}$, where $o\in \mathbf{T}_d$ is a root vertex.

	\subsection{Eigenvalue distribution}\label{s5.2}

	We now state and prove our eigenvalue distribution result, which by Lemma~\ref{lem:unitary-equiv} is equivalent to Theorem~\ref{t1.2}.
	\begin{theo}\label{t5.1}
		In the setting of~\eqref{n2.6}, assume \textbf{[BIV]}~\eqref{5.10}. Then for every interval $I\subseteq\mathbb{R}$,
		\begin{equation}\label{4.29}
			\lim_{\substack{N \in 2\N \\
					N \to\infty}}\lv\frac{1}{\dim_\mathbb{C}F_N}\lv\{i\mid\lambda_{N,i}\in I\}\rv-\int_{I}d\mu_{\mathrm{KM}}(\lambda)\rv=0.
		\end{equation}
	\end{theo}
	
	\begin{pro}
		By~\eqref{n2.9} and~\eqref{4.24}, we obtain	for $k\in\mathbb{N}$,
		\begin{equation}\label{4.26'}
			\int \lambda^kd\mu_{N}=\frac{1}{\dim_\mathbb{C}F_N}\sum_{{v}\in X,\gamma\in \Omega(X,{v}),\lv\gamma\rv=k}\mathrm{Tr}^{F_{N,{v}}}\big[\rho_N(\gamma)\big],
		\end{equation}
		
		We split the sum in~\eqref{4.26'} according to whether the homotopy class $[\gamma]\in\Gamma$, defined via the quotient in~\eqref{x2.4}, is trivial or nontrivial.
		\begin{equation}\label{4.17}
			\begin{split}
				\sum_{\substack{v\in X,\gamma\in \Omega(X,{v}),\\ \lv\gamma\rv=k}}\mathrm{Tr}^{F_{N,{v}}}\big[\rho_N(\gamma)\big]=&\sum_{\substack{v\in X,\gamma\in \Omega(X,{v}),\\
						\lv\gamma\rv=k,[\gamma]=1_\Gamma}}\mathrm{Tr}^{F_{N,{v}}}\big[\rho_N(\gamma)\big]\\
				&+\sum_{\substack{{v}\in X,\gamma\in \Omega(X,{v}), \\
						\lv\gamma\rv=k, [\gamma]\neq1_\Gamma }}\mathrm{Tr}^{F_{N,{v}}}\big[\rho_N(\gamma)\big].
			\end{split}
		\end{equation}
		If $[\gamma]=1_\Gamma$, the identity element of $\Gamma$, that is, the loop $\gamma$ is homotopic to a trivial loop, then $\gamma$ corresponds to a loop from $o$ to $o$ in $\mathbf{T}_d$. Therefore,
		\begin{equation}\label{4.18}
			\begin{split}
				\sum_{\substack{\gamma\in \Omega(X,{v}),\\
						\lv\gamma\rv=k,[\gamma]= 1_\Gamma}}\mathrm{Tr}^{F_{N,{v}}}\big[\rho_N(\gamma)\big]&=\sum_{\substack{\gamma\in \Omega(X,{v}),\\
						\lv\gamma\rv=k,[\gamma]= 1_\Gamma}}\mathrm{Tr}^{F_{N,{v}}}\big[\mathrm{Id}_{\mathcal{H}_N}\big]\\
				&=\bv\{\gamma\in\Omega(\mathbf{T}_d,o)\mid\lv\gamma\rv=k\}\bv\cdot \dim_{\mathbb{C}}\mathcal{H}_{N}.
			\end{split}
		\end{equation}
		If $[\gamma]\neq1_\Gamma$, then by \textbf{[BIV]}~\eqref{5.10}, we get from~\eqref{3.28x} and~\eqref{4.2x} that
		\begin{equation}\label{4.19}
			\begin{split}
				\lim_{N\to \infty}\bbv\frac{1}{\dim_{\mathbb{C}}\mathcal{H}_{N}}\mathrm{Tr}^{F_{N,{v}}}\big[\rho_N(\gamma)\big]\bbv&=\lim_{N\to \infty}\bbv\frac{1}{\dim_{\mathbb{C}}\mathcal{H}_{N}}\mathrm{Tr}^{F_{N,{v}}}\big[M_{\chi(\gamma),N}\big]\bbv\\
				&=0.
			\end{split}
		\end{equation}
		In other words, by the first pseudolocality result ~\eqref{3.28x}, the nontrivial loops do not contribute in the limit.

		By~\eqref{4.17},~\eqref{4.18}, and~\eqref{4.19}, we get for $k\in\mathbb{N}$,
		\begin{equation}
			\lim_{N\to \infty}\bbv\int \lambda^kd\mu_{N}(\lambda)-\bv\{\gamma\in \Omega(\mathbf{T}_d,o)\mid\lv\gamma\rv=k\}\bv\bbv=0,
		\end{equation}	
		and~\eqref{4.29} follows by~\eqref{eq:moment-identity} and an approximation argument.\qed
	\end{pro}

	\section{Kernel Operators}\label{S6}
	
	In this section, we introduce kernel operators associated with unitary vector bundles over graphs. In \S\,\ref{s7.1nn}, we introduce the endomorphism vector bundle. In \S\,\ref{s7.1n}, we enlarge it to the kernel vector bundle. In \S\,\ref{s6.3x}, we establish several basic properties of kernel operators. In \S\,\ref{s7.3n}, we study the commutator and gradient operators.
	
	In this section, we mainly refer to ~\cite[\S\S\,2, 3, 5]{MR3649482} and ~\cite[\S\,8]{disc}.
	\subsection{Endomorphism vector bundles}\label{s7.1nn}
	
	Choose $(X,F)$, where $X$ is a $d$-regular graph and $F$ is a unitary flat vector bundle of the form~\eqref{n2.4}. 
	
	Define the endomorphism bundle $\mathrm{End}(F)$ by
	\begin{equation}\label{7.58}
		\begin{split}
			\mathrm{End}(F)& \coloneqq \Gamma\backslash\big(\mathbf{T}_d\times\mathrm{End}(\mathbb{C}^\ell)\big)\\
			&=\big\{(\bm{v},A)\in \mathbf{T}_d\times\mathrm{End}(\mathbb{C}^\ell) \big\}\big/(\bm{v},A)\sim(\gamma\bm{v},\rho(\gamma)A\rho(\gamma)^{-1})\\
			&\qquad\qquad\qquad\qquad\qquad\qquad\quad\text{for every }\gamma\in\Gamma.
		\end{split}
	\end{equation}
	Then, as in~\eqref{n2.5}, we have
	\begin{equation}\label{7.59}
		\begin{split}
			C\big(X, \mathrm{End}(F)\big)=\big\{Q\colon \mathbf{T}_d\to \mathrm{End}(\mathbb{C}^\ell)\mid\ &Q(\gamma \bm{v})=\rho(\gamma)Q(\bm{v})\rho(\gamma)^{-1}\\
			&\text{for every } \bm{v}\in \mathbf{T}_d, \gamma\in\Gamma\big\}.
		\end{split}
	\end{equation}
	
	The space $C\big(X, \mathrm{End}(F)\big)$ is too small for the analysis below, in particular, it is not stable under composition with $\Delta^F$. We therefore enlarge the class of operators under consideration.

	\subsection{Kernel vector bundles}\label{s7.1n}

	Define \emph{the kernel space $K$} and \emph{the kernel vector bundle $K^{\mathrm{End}(F)}$} by
	\begin{equation}\label{7.1n}
		\begin{split}
			K&\coloneqq\Gamma\backslash \mathbf{T}_d^{\times 2}\\
			&=\big\{(\bm{v},\bm{v}')\in \mathbf{T}_d^{\times 2}\big\}\big/(\bm{v},\bm{v}')\sim(\gamma \bm{v},\gamma \bm{v}')\ \text{for every }\gamma\in\Gamma,\\
			K^{\mathrm{End}(F)}& \coloneqq \Gamma\backslash\big(\mathbf{T}_d^{\times 2}\times \mathrm{End}(\mathbb{C}^\ell)\big)\\
			&=\big\{(\bm{v},\bm{v}',A)\in \mathbf{T}_d^{\times 2}\times\mathrm{End}(\mathbb{C}^\ell)\big\}\big/(\bm{v},\bm{v}',A)\sim(\gamma \bm{v},\gamma \bm{v}',\rho(\gamma)A\rho(\gamma)^{-1})\\
			&\qquad\qquad\qquad\qquad\qquad\qquad\qquad\quad\text{for every }\gamma\in\Gamma.
		\end{split}
	\end{equation}
	We call a section $Q(\bm{v},\bm{v}')\in C(K, K^{\mathrm{End}(F)})$ a \emph{kernel operator}. As in~\eqref{n2.5}, we can identify $C(K, K^{\mathrm{End}(F)})$ with the space of $\Gamma$-equivariant $\mathrm{End}(\mathbb{C}^\ell)$-valued functions on $\mathbf{T}_d^{\times 2}$, namely
	\begin{equation}\label{7.2n}
		\begin{split}
			C(K, K^{\mathrm{End}(F)})=\big\{Q\colon\mathbf{T}_d^{\times 2}\to \mathrm{End}(\mathbb{C}^\ell) \mid\ &Q(\gamma \bm{v},\gamma  \bm{v}')=\rho(\gamma)Q(\bm{v},\bm{v}')\rho(\gamma)^{-1}\\
			&\text{for every } (\bm{v},\bm{v}')\in \mathbf{T}_d^{\times 2}, \gamma\in\Gamma\big\}.
		\end{split}
	\end{equation}

	Notice that a nonbacktracking path $(\bm{v}_k,\cdots,\bm{v}_1,\bm{v}_0)$ in $\mathbf{T}_d$ is uniquely determined by its endpoints $(\bm{v}_k,\bm{v}_0)\in \mathbf{T}_d^{\times 2}$. From this viewpoint, we may regard $K$, defined in~\eqref{7.1n}, as the space of nonbacktracking paths in $X$. We denote by $K_{(k)}$ the space of paths of length $k$, and by $K_{(\leqslant k)}$ the space of paths of length at most $k$,
	\begin{equation}\label{7.4}
		\begin{split}
			K_{(k)}& \coloneqq \{\text{nonbacktracking paths of length $k$ in }X\}\\
			&\cong\{(\bm{v},\bm{v}')\in K\mid d_{\mathbf{T}_d}(\bm{v},\bm{v}')=k\},\\
			K_{(\leqslant k)}& \coloneqq \mathsmaller{\bigcup}_{i=0}^{k}K_{(i)}.
		\end{split}
	\end{equation}
	Accordingly, we have the direct sum decomposition of~\eqref{7.2n},
	\begin{equation}\label{7.5nnn}
		\begin{split}
			C(K, K^{\mathrm{End}(F)})&=\mathsmaller{\bigoplus}_{k=0}^\infty C(K_{(k)},K^{\mathrm{End}(F)}),\\ Q&=Q_{(0)}+Q_{(1)}+\cdots,
		\end{split}
	\end{equation}
	where for $Q\in C(K, K^{\mathrm{End}(F)})$, we denote by $Q_{(k)}\in C(K_{(k)}, K^{\mathrm{End}(F)})$ the component of $Q$ supported on $K_{(k)}$. In particular, by~\eqref{7.59}, we have
	\begin{equation}\label{7.10nn}
		C(K_{(0)}, K^{\mathrm{End}(F)})\cong C(X,\mathrm{End}(F)).
	\end{equation}
	Hence $C(K,K^{\mathrm{End}(F)})$ provides the desired enlargement of $C(X,\mathrm{End}(F))$.

	\subsection{Properties of kernel operators}\label{s6.3x}

	By~\eqref{n2.5} and~\eqref{7.2n}, a kernel operator $Q\in C(K, K^{\mathrm{End}(F)})$ acts on a section $u\in C(X,F)$ by
	\begin{equation}\label{7.5n}
		\begin{split}
			(Qu)(\bm{v})=&\sum_{\bm{v}'\in \mathbf{T}_d}Q(\bm{v},\bm{v}')u(\bm{v}')=\sum_{\bm{v}'\in D,\gamma\in\Gamma}Q(\bm{v},\gamma \bm{v}')u(\gamma \bm{v}')\\
			=&\sum_{\bm{v}'\in D,\gamma\in\Gamma}Q(\bm{v},\gamma \bm{v}')\rho(\gamma)u(\bm{v}'),
		\end{split}
	\end{equation}
	where $D$ is a vertex fundamental domain as in \eqref{n2.1}.

	\begin{lemma}
		For two kernel operators $Q,Q'\in C(K,K^{\mathrm{End}(F)})$, their composition $QQ'$ is again a kernel operator with
		\begin{equation}\label{7.6n}	(QQ')(\bm{v},\bm{v}')=\sum_{\bm{v}''\in\mathbf{T}_d}Q(\bm{v},\bm{v}'')Q'(\bm{v}'',\bm{v}').
		\end{equation}
	\end{lemma}
	
	\begin{pro}
		By~\eqref{7.5n}, we have
		\begin{equation}
			\begin{split}
				\big(QQ'u\big)(\bm{v})&=\sum_{\bm{v}'\in\mathbf{T}_d}Q(\bm{v},\bm{v}')(Q'u)(\bm{v}')\\
				&=\sum_{\bm{v}',\bm{v}''\in\mathbf{T}_d}Q(\bm{v},\bm{v}')Q'(\bm{v}',\bm{v}'')u(\bm{v}''),
			\end{split}
		\end{equation}
		which implies~\eqref{7.6n}.\qed
	\end{pro}

	By~\eqref{eq:hermitian -induces} and~\eqref{7.5n}, we have
	\begin{equation}\label{7.14nn}
		\begin{split}
			\big\langle Qu,u'\big\rangle_{L^2(X,F)}=&\sum_{\bm{v}\in D}\big\langle(Qu)(\bm{v}),u'(\bm{v})\big\rangle_{\mathbb{C}^\ell}\\
			=&\sum_{\bm{v},\bm{v}'\in D,\gamma\in\Gamma}\big\langle Q(\bm{v},\gamma \bm{v}')\rho(\gamma)u(\bm{v}'), u'(\bm{v})\big\rangle_{\mathbb{C}^\ell}.
		\end{split}
	\end{equation}

	We denote by $Q^*$ the adjoint of $Q$ on $L^2(X,F)$. We define the normalized
	Hilbert--Schmidt norm of $Q$ by
	\begin{equation}\label{6.12}
		\lV Q\rV_{\mathrm{HS}(X,F)}^2 \coloneqq \frac{1}{\dim_{\mathbb{C}}F}\mathrm{Tr}^{L^2(X,F)}\big[Q^*Q\big],
	\end{equation}
	where we use the dimension conventions in~\eqref{x2.5}.
	
	\begin{lemma}
		For $Q\in C(K,K^{\mathrm{End}(F)})$, we have
		\begin{equation}\label{7.9n}
			\lV Q\rV_{\mathrm{HS}(X,F)}^2=\frac{1}{\dim_{\mathbb{C}}F}\sum_{\bm{v},\bm{v}'\in D}\BV\sum_{\gamma\in\Gamma}Q(\bm{v},\gamma \bm{v}')\rho(\gamma)\BV_{\mathrm{HS}(\mathbb{C}^\ell)}^2,
		\end{equation}
		where the Hilbert--Schmidt norm on
		$\mathrm{End}(\mathbb{C}^\ell)$ is defined by
		\begin{equation}
			\lV Q(\bm{v},\bm{v}')\rV_{\mathrm{HS}(\mathbb{C}^\ell)}^2=\mathrm{Tr}^{\mathbb{C}^\ell}\big[Q(\bm{v},\bm{v}')^*\cdot Q(\bm{v},\bm{v}')\big].
		\end{equation}
	\end{lemma}
	
	\begin{pro}
		Choose an orthonormal basis of $L^2(X,F)$ consisting of sections supported at individual points of $D$, with values in an orthonormal basis of $\mathbb{C}^\ell$. By~\eqref{7.14nn}, the matrix block of $Q$ indexed by $\bm{v},\bm{v}'\in D$ is $\sum_{\gamma\in\Gamma}Q(\bm{v},\gamma \bm{v}')\rho(\gamma)$. Summing the squared Hilbert--Schmidt norms of these blocks and dividing by the dimension gives~\eqref{7.9n}.\qed
	\end{pro}

	As we will see, quantum ergodicity requires control of $\lV\cdot\rV_{\mathrm{HS}(X,F)}$, which is difficult to compute directly, as it involves correlations between distinct points. Thus, we introduce another norm that is easier to compute. We define the normalized $L^{2}$-norm of $Q$ by
	\begin{equation}\label{7.7n}
		\begin{split}
			\lV Q\rV_{L^2(K,K^{\mathrm{End}(F)})}^2 \coloneqq \frac{1}{\dim_{\mathbb{C}}F}\sum_{\bm{v}\in D,\bm{v}'\in \mathbf{T}_d}\lV Q(\bm{v},\bm{v}')\rV_{\mathrm{HS}(\mathbb{C}^\ell)}^2.
		\end{split}
	\end{equation}

	\begin{lemma}
		Let $Q\in C(K, K^{\mathrm{End}(F)})$. Then
		\begin{equation}\label{7.10n}
			\begin{split}
				\lV Q\rV_{\mathrm{HS}(X,F)}^2=&\lV Q\rV_{L^2(K,K^{\mathrm{End}(F)})}^2\\
				&+\sum_{\substack{\bm{v},\bm{v}'\in D,\  \gamma,\gamma'\in\Gamma\\ \gamma\neq \gamma'}}\frac{\mathrm{Tr}^{\mathbb{C}^\ell} }{\dim_{\mathbb{C}}F}\Big[\rho(\gamma')^{-1}Q(\bm{v},\gamma'\bm{v}')^*Q(\bm{v},\gamma \bm{v}')\rho(\gamma)\Big].
			\end{split}
		\end{equation}
	\end{lemma}
	
	\begin{pro}
		By direct computation, for $(\bm{v},\bm{v}')\in \mathbf{T}_d^{\times 2}$,
		\begin{equation}\label{7.11n}
			\begin{split}
				&\BV\sum_{\gamma\in\Gamma}Q(\bm{v},\gamma \bm{v}')\rho(\gamma)\BV_{\mathrm{HS}(\mathbb{C}^\ell)}^2\\
				&=\sum_{\gamma\in\Gamma}\BV Q(\bm{v},\gamma \bm{v}')\BV_{\mathrm{HS}(\mathbb{C}^\ell)}^2+\sum_{\gamma\neq \gamma'}\mathrm{Tr}^{\mathbb{C}^\ell}\Big[\rho(\gamma')^{-1}Q(\bm{v},\gamma'\bm{v}')^*Q(\bm{v},\gamma \bm{v}')\rho(\gamma)\Big].
			\end{split}
		\end{equation}
		Summing~\eqref{7.11n} over $\bm{v},\bm{v}'\in D$, by~\eqref{7.9n} and~\eqref{7.7n}, we get~\eqref{7.10n}.\qed
	\end{pro}

	Our overall strategy is to bound $\lV\cdot\rV_{\mathrm{HS}(X,F)}$ in terms of $\lV\cdot\rV_{L^2(K,K^{\mathrm{End}(F)})}$. We shall formulate assumptions tailored to ensure that the correlation term on the right hand side of~\eqref{7.10n} is asymptotically negligible.

	\subsection{Commutator and gradient operators}\label{s7.3n}
	
	We now discuss several operators acting on $C(K,K^{\mathrm{End}(F)})$. Throughout this subsection, adjoints are taken with respect to the inner product associated with the norm $\lV\cdot\rV_{L^2(K,K^{\mathrm{End}(F)})}$ defined in~\eqref{7.7n}.

	Let $Q\in C(K,K^{\mathrm{End}(F)})$ be a kernel operator.  Its commutator $[\Delta^F,Q]$ with $\Delta^F$ is again a kernel operator. This commutator will play an important role in our discussion of quantum ergodicity.
	
	Using~\eqref{7.6n}, we define the commutator map
	\begin{equation}\label{7.27n}
		\begin{split}
			&\mathrm{ad}_{\Delta^F}\colon C(K_{(k)},K^{\mathrm{End}(F)}) \to C(K_{(k-1)}\cup K_{(k+1)},K^{\mathrm{End}(F)}),\\
			&\mathrm{ad}_{\Delta^F}(Q)(\bm{v},\bm{v}')=\big({\Delta^F}Q-Q{\Delta^F}\big)(\bm{v},\bm{v}')\\
			&\qquad\qquad\qquad\ =\sum_{\bm{v}''\sim \bm{v}}Q(\bm{v}'',\bm{v}')-\sum_{\bm{v}'''\sim \bm{v}'}Q(\bm{v},\bm{v}''').
		\end{split}
	\end{equation}
	Indeed, if $Q\in C(K_{(k)},K^{\mathrm{End}(F)})$, then
	\begin{equation}\label{7.22n}
		\begin{split}
			&(\Delta^FQ)(\bm{v}_{k-1},\cdots,\bm{v}_0)=\sum_{\bm{v}_k\sim \bm{v}_{k-1}, \bm{v}_k\neq \bm{v}_{k-2}}Q(\bm{v}_{k},\cdots,\bm{v}_0),\\
			&(\Delta^FQ)(\bm{v}_{k+1},\bm{v}_k,\cdots,\bm{v}_0)=Q(\bm{v}_{k},\cdots,\bm{v}_0),\\
			&(Q\Delta^F)(\bm{v}_k,\cdots,\bm{v}_1)=\sum_{\bm{v}_0\sim \bm{v}_{1}, \bm{v}_0\neq \bm{v}_2}Q(\bm{v}_k,\cdots,\bm{v}_0),\\
			&(Q\Delta^F)(\bm{v}_{k},\cdots,\bm{v}_0,\bm{v}_{-1})=Q(\bm{v}_{k},\cdots,\bm{v}_0).
		\end{split}
	\end{equation}

	By~\eqref{7.27n} and~\eqref{7.22n}, we can decompose $\mathrm{ad}_{\Delta^F}$ into two components,
	\begin{equation}\label{7.33nn}
		\mathrm{ad}_{\Delta^F}=-(\nabla+\nabla^*).
	\end{equation}
	Here $\nabla$ is the gradient operator
	\begin{equation}\label{7.33n}
		\begin{split}
			&\nabla\colon C(K_{(k)},K^{\mathrm{End}(F)}) \to C(K_{(k+1)},K^{\mathrm{End}(F)}),\\
			&\nabla(Q)(\bm{v}_{k+1},\cdots,\bm{v}_{0})=Q(\bm{v}_{k+1},\cdots,\bm{v}_{1})-Q(\bm{v}_k,\cdots,\bm{v}_0),
		\end{split}
	\end{equation}
	and $\nabla^*$ is its adjoint (as we will verify in Lemma~\ref{lem:nabla-adjoint}),
	\begin{equation}\label{7.34n}
		\begin{split}
			&\nabla^*\colon C(K_{(k)},K^{\mathrm{End}(F)}) \to C(K_{(k-1)},K^{\mathrm{End}(F)}),\\
			&\nabla^*(Q)(\bm{v}_k,\cdots,\bm{v}_1)=\sum_{\bm{v}_0\sim \bm{v}_1,\bm{v}_0\neq \bm{v}_2}Q(\bm{v}_k,\cdots,\bm{v}_{0})\\
			&\qquad\qquad\qquad\qquad\quad-\sum_{\bm{v}_{k+1}\sim \bm{v}_k,\bm{v}_{k+1}\neq \bm{v}_{k-1}}Q(\bm{v}_{k+1},\cdots,\bm{v}_{1}).
		\end{split}
	\end{equation}
	Note that when $k=1$, formula~\eqref{7.34n} becomes
	\begin{equation}\label{7.34nn}
		(\nabla^*Q)(\bm v_1)=\sum_{\bm v_0\sim\bm v_1}Q(\bm v_1,\bm v_0)-\sum_{\bm v_2\sim\bm v_1}Q(\bm v_2,\bm v_1).
	\end{equation}
	
	We now validate the definition of $\nabla^*$.
	\begin{lemma}\label{lem:nabla-adjoint}
	On $L^2(K,K^{\mathrm{End}(F)})$, the operator $\nabla^*$ is the adjoint of the operator $\nabla$. 
	\end{lemma}
	
	\begin{pro}
		
		Let $Q_{(k-1)}\in L^2(K_{(k-1)},K^{\mathrm{End}(F)})$ and $Q_{(k)}\in L^2(K_{(k)},K^{\mathrm{End}(F)})$. By~\eqref{7.7n} and~\eqref{7.33n}, we can split $\big\langle\nabla Q_{(k-1)},Q_{(k)}\big\rangle_{L^2(K,K^{\mathrm{End}(F)})}$ into the difference of the following two terms,
		\begin{equation}\label{6.24}
			\begin{split}
				&\langle \nabla Q_{(k-1)}, Q_{(k)} \rangle_{L^2(K,K^{\mathrm{End}(F)})} \\
				&= \frac{1}{\dim_{\mathbb C}F}
				\sum_{\substack{\bm v_k\in D,\\
						(\bm v_k,\cdots,\bm v_0)}}
				\big\langle
				Q_{(k-1)}(\bm v_k,\cdots,\bm v_1),
				Q_{(k)}(\bm v_k,\cdots,\bm v_0)
				\big\rangle_{\mathrm{HS}(\mathbb C^\ell)},\\
				& \quad -\frac{1}{\dim_{\mathbb C}F}
				\sum_{\substack{\bm v_k\in D,\\
						(\bm v_k,\cdots,\bm v_0)}}\big\langle
				Q_{(k-1)}(\bm v_{k-1},\cdots,\bm v_0),
				Q_{(k)}(\bm v_k,\cdots,\bm v_0)
				\big\rangle_{\mathrm{HS}(\mathbb C^\ell)},
			\end{split}
		\end{equation}
		where the sum is over nonbacktracking paths $(\bm v_k,\cdots,\bm v_0)$.
		
		For the first term in the RHS of~\eqref{6.24}, we fix the shortened path $(\bm v_k,\cdots,\bm v_1)$ and then sum over all its nonbacktracking extensions to get
		\begin{equation}\label{7.35nb}
			\begin{split}
				&\frac{1}{\dim_{\mathbb C}F}
				\sum_{\substack{\bm v_{k}\in D,\\
						(\bm v_{k},\cdots,\bm v_1)}}
				\Big\langle Q_{(k-1)}(\bm v_{k},\cdots,\bm v_1),\sum_{\substack{\bm v_{0}\sim\bm v_1,\\
						\bm v_{0}\neq\bm v_2}}
				Q_{(k)}(\bm v_{k},\cdots,\bm v_0)\Big\rangle_{\mathrm{HS}(\mathbb C^\ell)}.
			\end{split}
		\end{equation}
		
		For the second term in the RHS of~\eqref{6.24}, the vertex $\bm v_{k-1}$ need not belong to $D$. The diagonal $\Gamma$-action gives a bijection
		\begin{equation}\label{7.35nd}
			\big\{(\bm v_{k},\cdots,\bm v_0)\mid\bm v_{k}\in D\big\}\cong\big\{(\bm v_{k},\cdots,\bm v_0)\mid\bm v_{k-1}\in D\big\}.
		\end{equation}
		Indeed, given a path
		$(\bm v_{k},\cdots,\bm v_0)$ with $\bm v_{k}\in D$, there exists a unique $\gamma\in\Gamma$ such that $\gamma\bm v_{k-1}\in D$. We associate with it the path $(\gamma\bm v_{k},\cdots,\gamma\bm v_0)$. Conversely, starting with a path with $\bm v_{k-1}\in D$, there exists a unique $\gamma\in \Gamma$ such that $\gamma\bm v_{k}\in D$. 
		
		By the $\Gamma$-equivariance of $Q_{(k-1)}$ and $Q_{(k)}$ in~\eqref{7.2n}, and the unitarity of $\rho$ in~\eqref{unitary}, we get
		\begin{equation}\label{7.35nd'}
			\begin{split}
				&\big\langle
				Q_{(k-1)}(\gamma\bm v_{k-1},\cdots,\gamma\bm v_0),
				Q_{(k)}(\gamma\bm v_k,\cdots,\gamma\bm v_0)
				\big\rangle_{\mathrm{HS}(\mathbb C^\ell)}\\
				&=\big\langle \rho(\gamma)Q_{(k-1)}(\bm v_{k-1},\cdots,\bm v_0)\rho(\gamma)^{-1},
				\rho(\gamma)Q_{(k)}(\bm v_k,\cdots,\bm v_0)\rho(\gamma)^{-1}
				\big\rangle_{\mathrm{HS}(\mathbb C^\ell)}\\
				&=\big\langle
				Q_{(k-1)}(\bm v_{k-1},\cdots,\bm v_0),
				Q_{(k)}(\bm v_k,\cdots,\bm v_0)
				\big\rangle_{\mathrm{HS}(\mathbb C^\ell)}.
			\end{split}
		\end{equation}
		Consequently, the second term in~\eqref{6.24} is equal to
		\begin{equation}\label{7.35nc}
			\begin{split}
				\frac{1}{\dim_{\mathbb C}F}\sum_{\substack{\bm v_{k-1}\in D,\\
						(\bm v_{k-1},\cdots,\bm v_0)}}\Big\langle
				Q_{(k-1)}(\bm v_{k-1},\cdots,\bm v_0),
				\sum_{\substack{\bm v_k\sim\bm v_{k-1},\\ \bm v_k\neq\bm v_{k-2}}}
				Q_{(k)}(\bm v_k,\cdots,\bm v_0)
				\Big\rangle_{\mathrm{HS}(\mathbb C^\ell)}.
			\end{split}
		\end{equation}
		
		Combining~\eqref{6.24},~\eqref{7.35nb}, and~\eqref{7.35nc}, we finish the proof.\qed
	\end{pro}

	\begin{lemma}
		The commutator is selfadjoint acting on $L^2(K, K^{\mathrm{End}(F)})$,
		\begin{equation}\label{7.36}
			\mathrm{ad}_{\Delta^F}^*=\mathrm{ad}_{\Delta^F}.
		\end{equation}
	\end{lemma}
	
	\begin{pro}
		This follows immediately from~\eqref{7.33nn}.\qed
	\end{pro}

	We define a kernel operator $\mathrm{Id}_{(k)}$ for every $k\in\mathbb{N}$ by
	\begin{equation}\label{7.5nn}
		\mathrm{Id}_{(k)}(\bm{v},\bm{v}')=\mathrm{Id}_{\mathbb{C}^\ell}\mathbbm{1}_{\{d_{\mathbf{T}_d}(\bm{v},\bm{v}')=k\}}(\bm{v},\bm{v}').
	\end{equation}

	\begin{prop}\label{L7.17}
On $L^2(K,K^{\mathrm{End}(F)})$, operators $\mathrm{ad}_{\Delta^F}$ and $\nabla$ satisfy
		\begin{equation}\label{7.87}
			\begin{split}
				\ker(\mathrm{ad}_{\Delta^F})&=\ker(\nabla)\\
				&=\mathrm{span}\big\{A\cdot\mathrm{Id}_{(k)}\mid k\in\mathbb{N},  A\in\mathrm{End}(\mathbb{C}^\ell),\\
				&\qquad\qquad\qquad\qquad\rho(\gamma)A\rho(\gamma)^{-1}=A \textnormal{ for every }\gamma\in \Gamma\big\},
			\end{split}
		\end{equation}
	where the span is understood as the closed linear span.
	\end{prop}
	
	\begin{pro}
		
		We first claim that for $k\geqslant 2$,
		\begin{equation}\label{7.38n}
			\bV\nabla^*(Q_{(k)})\bV_{L^2(K_{(k-1)},K^{\mathrm{End}(F)})}^2\leqslant\bV\nabla (Q_{(k)})\bV_{L^2(K_{(k+1)},K^{\mathrm{End}(F)})}^2.
		\end{equation}
		By~\eqref{7.34n}, we have
		\begin{equation}\label{6.33}
			\begin{split}
				&(\nabla^* Q_{(k)})(\bm{v}_k,\cdots,\bm{v}_1)\\
				&=\frac{1}{d-1}\sum_{\substack{\bm{v}_0\sim \bm{v}_1,\bm{v}_0\neq \bm{v}_2\\ \bm{v}_{k+1}\sim \bm{v}_k,\bm{v}_{k+1}\neq \bm{v}_{k-1}
				}}\big(Q_{(k)}(\bm{v}_k,\cdots,\bm{v}_0)-Q_{(k)}(\bm{v}_{k+1},\cdots,\bm{v}_1)\big).
			\end{split}
		\end{equation}
		Here the factor $\frac{1}{d-1}$ appears because, in \eqref{7.34n}, both the sum over $\bm{v}_0$ and the sum over $\bm{v}_{k+1}$ contain $(d-1)$ terms. When these two separate sums are rewritten as a single double sum, each summand is counted $(d-1)$ times.	By ~\eqref{7.33n}, the right hand side of \eqref{6.33} can be written as
		\begin{equation}
			\frac{1}{d-1}\sum_{\substack{\bm{v}_0\sim \bm{v}_1,\bm{v}_0\neq \bm{v}_2\\ \bm{v}_{k+1}\sim \bm{v}_k,\bm{v}_{k+1}\neq \bm{v}_{k-1}
			}}-\nabla (Q_{(k)})(\bm{v}_{k+1},\cdots,\bm{v}_{0}).
		\end{equation}
		
		Applying the Cauchy--Schwarz inequality, we get
		\begin{equation}
			\begin{split}
				&\lV(\nabla^* Q_{(k)})(\bm{v}_k,\cdots,\bm{v}_1)\rV_{\mathrm{HS}(\mathbb{C}^\ell)}^2\\
				&\leqslant\sum_{\substack{\bm{v}_0\sim \bm{v}_1,\bm{v}_0\neq \bm{v}_2\\ \bm{v}_{k+1}\sim \bm{v}_k,\bm{v}_{k+1}\neq \bm{v}_{k-1}
				}}\lV\nabla (Q_{(k)})(\bm{v}_{k+1},\cdots,\bm{v}_{0})\rV_{\mathrm{HS}(\mathbb{C}^\ell)}^2.
			\end{split}
		\end{equation}
		Then~\eqref{7.38n} follows from~\eqref{7.7n} and an argument similar to that in~\eqref{7.35nd} and~\eqref{7.35nd'}.

		By~\eqref{7.33nn}, if $\mathrm{ad}_{\Delta^F}Q=0$, taking the $L^2(K_{(k+1)},K^{\mathrm{End}(F)})$ component ($k \geqslant 0$) gives
		\begin{equation}\label{6.35}
			\nabla \big(Q_{(k)}\big)+\nabla^* \big(Q_{(k+2)}\big)=0.
		\end{equation}
		From~\eqref{7.38n},~\eqref{6.35}, and the fact that $(k+2)\geqslant 2$, we get
		\begin{equation}\label{6.36}
			\begin{split}
				\bV\nabla(Q_{(k)})\bV_{L^2(K_{(k+1)},K^{\mathrm{End}(F)})}^2&=\bV\nabla^* \big(Q_{(k+2)}\big)\bV_{L^2(K_{(k+1)},K^{\mathrm{End}(F)})}^2\\		
				&\leqslant\bV\nabla  (Q_{(k+2)})\bV_{L^2(K_{(k+3)},K^{\mathrm{End}(F)})}^2.
			\end{split}
		\end{equation}
		Since $Q\in L^2(K,K^{\mathrm{End}(F)})$ implies that $\nabla Q\in L^2(K,K^{\mathrm{End}(F)})$, we have
		\begin{equation}\label{6.37x}
			\lim_{k\to\infty}\Vert\nabla (Q_{(k)})\Vert_{L^2(K_{(k+1)},K^{\mathrm{End}(F)})}^2=0.
		\end{equation}
		By \eqref{6.36} and \eqref{6.37x}, we have $\nabla Q=0$ and $\ker(\mathrm{ad}_{\Delta^F})\subseteq\ker(\nabla)$. By~\eqref{7.33n}, if $\nabla Q_{(k)}=0$, then $Q_{(k)}$ is actually a constant matrix, taking into account the $\Gamma$-equivariance in~\eqref{7.2n}, we get the second identity of~\eqref{7.87}. For \(Q_{(k)}=A\cdot\mathrm{Id}_{(k)}\), clearly $\nabla Q_{(k)}=0$ and $\nabla^*Q_{(k)}=0$ and $\ker(\nabla)\subseteq \ker(\mathrm{ad}_{\Delta^F})$. We get the first identity of~\eqref{7.87}.\qed
	\end{pro}

	Finally, we show that the kernel operator $\mathrm{Id}_{(k)}$ defined in~\eqref{7.5nn} can be written in terms of $\Delta^F$. We define a sequence of polynomials $(h_k(\lambda))_{k\in\mathbb{N}}$ by
	\begin{equation}\label{6.37}
		\begin{split}
			&h_0(\lambda)=1,\qquad
			h_1(\lambda)=\lambda,\qquad
			h_2(\lambda)=\lambda^2-d,\\
			&h_{k+1}(\lambda)=\lambda h_k(\lambda)-(d-1)h_{k-1}(\lambda), \text{ for } k\geqslant2.
		\end{split}
	\end{equation}
	
	\begin{prop}
		We have for $k\in\mathbb{N}$,
		\begin{equation}\label{6.31}
			h_k(\Delta^F)=\mathrm{Id}_{(k)}.
		\end{equation}
	\end{prop}
	\begin{pro}
		
		Clearly we have
		\begin{equation}\label{6.39}
			\mathrm{Id}_{(0)}=\mathrm{Id}_F,
			\qquad
			\mathrm{Id}_{(1)}=\Delta^F.
		\end{equation}
		Moreover,
		\begin{equation}\label{6.40}
			\begin{split}
				&\Delta^F\mathrm{Id}_{(1)}=\mathrm{Id}_{(2)}+d\mathrm{Id}_{(0)},\\
				&\Delta^F\mathrm{Id}_{(k)}=\mathrm{Id}_{(k+1)}+(d-1)\mathrm{Id}_{(k-1)} \text{ for } k\geqslant2.
			\end{split}
		\end{equation}
		We see this from the following interpretation of the coefficients $d$ and $(d-1)$. Composing $\mathrm{Id}_{(k)}$ with $\Delta^F$ amounts to extending a nonbacktracking path of length $k$ by one edge. The resulting path either remains nonbacktracking and has length $k+1$, or its last two edges cancel. When $k=1$, there are $d$ paths that cancel to a single vertex. When $k\geqslant2$, each nonbacktracking path of length $(k-1)$ arises from $(d-1)$ such cancellations.
		
		Combining~\eqref{6.39},~\eqref{6.40}, and ~\eqref{6.37}, we  conclude ~\eqref{6.31}.\qed
	\end{pro}

	\section{Quantum Cat Maps, Continued}\label{S7}
	
	In this section, we continue the discussion of \S\,\ref{s3}, where quantum cat maps were introduced as quantum evolution operators. We now introduce quantized functions on the torus, which play the role of quantum observables. In \S\,\ref{s7.1}, we recall Weyl quantization on the Euclidean space. In \S\,\ref{s7.2}, we introduce the quantization of functions on the torus. In \S\,\ref{s7.3}, we establish a pseudolocality estimate involving quantum cat maps and quantum observables. In \S\,\ref{s7.4}, we prove several trace asymptotics.

	\subsection{Quantization on the Euclidean space}\label{s7.1}
	
	We begin by recalling the semiclassical Weyl quantization. Let $\rS(\R^{n})$ be the Schwartz space of rapidly decreasing functions, with its dual space $\mathscr{S}'(\R^{n})$ of tempered distributions.

	For $a \in \rS(\R^{2n})$ and a semiclassical parameter $h \in (0,1]$, its Weyl quantization is defined by
	\begin{equation}\label{3.1f}
		\begin{split}
			&\op_h(a)\colon \rS(\R^n)\to \rS(\R^n),\\
			&\big(\op_h(a)f\big)(x)=\frac{1}{(2\pi h)^n} \int_{\R^{2n}} e^{\frac{\sqrt{-1}}{h} \langle x-{y},{\eta}\rangle}a\big(\tfrac{x+{y}}{2},{\eta}\big)f({y}) d{y} d{\eta}.
		\end{split}
	\end{equation}
	
	We then define the symbol class 
	\begin{equation}\label{3.2f}
		S(1) \coloneqq \Big\{a\in C^\infty\left(\R^{2n}\right)\mid\sup_{z\in\R^{2n}}\bv\partial^\alpha a(z)\bv<\infty \text{ for all } \alpha \in \N^{2n}\Big\}.
	\end{equation}
	From ~\cite[Theorem 4.16 and Theorem 4.23]{MR2952218}, when $a \in S(1)$, the operator $\op_h(a)$ acts on both $\mathscr{S}(\R^n)$ and  $\mathscr{S}'(\R^n)$ and is bounded on $L^2(\R^n)$. We now formally define the quantum translation from ~\eqref{3.1x}.  For $\omega \in\R^{2n}$, define
	\begin{equation}\label{3.3}
		U_\omega \coloneqq \op_h(a_\omega), \quad a_\omega(z) \coloneqq \exp(\tfrac{\sqrt{-1}}{h} \langle \omega,z\rangle), \quad z \in \R^{2n}.
	\end{equation}
	The explicit formula for the action of $U_\omega$ given in~\eqref{3.1x} follows from ~\cite[Theorem 4.7]{MR2952218}.

	Using~\eqref{3.1x}, we see that $U_\omega$ satisfies an exact Egorov theorem. For any $a\in S(1)$,
	\begin{equation}\label{3.5f}
		U_\omega^{-1}\op_h(a)U_\omega=\op_h(\tilde{a}), \quad \tilde{a}({x},{\xi})=a({x}-\eta,{\xi}+y).
	\end{equation}
	
	We now also formally define $\mathcal{M}_E$, the set of metaplectic transformations associated with $E \in \Sp_{2n}(\R)$. The set $\mathcal{M}_E$ is defined to be the set of all unitary transformations $M_E: L^2(\R^n) \rightarrow L^2(\R^n)$  that satisfy the following exact Egorov theorem, 
	\begin{equation}\label{MA}
		M_{E}^{-1} \op_h(a) M_E=\op_h(a \circ E).
	\end{equation}
	We once again remark that such transformations are unique up to multiplication by a unit complex scalar. 
	Clearly ~\eqref{3.8xx}, the exact Egorov theorem we used in the proof of Theorem~\ref{t1.2}, is a special case of~\eqref{MA}.

	\subsection{Quantization on the torus}\label{s7.2}

	We now assume that $a \in C^\infty(\bT^{2n})$. We can identify each such $a$ with a $\Z^{2n}$-periodic function on $\R^{2n}$. Since every $a \in C^\infty(\bT^{2n})$ is in the symbol class $S(1)$ defined in~\eqref{3.2f}, the quantization $\op_h(a)$ acts on $L^2(\R^n)$.

	From~\eqref{3.5f}, we have for any $\omega \in \Z^{2n}$,
	\begin{equation}\label{3.9f}
		\op_h(a)U_\omega=U_\omega \op_h(a).
	\end{equation}
	By~\eqref{3.9f}, we can check that under~\eqref{x3.9}, the operator $\op_h(a)$ acts on the space $\cH_N$ of quantum states defined in~\eqref{3.11f}. We denote this action by
	\begin{equation}\label{7.7}
		\op_N(a) \coloneqq \op_h(a)\big|_{\cH_N}.
	\end{equation} 
	
	From~\eqref{MA}, we see that
	\begin{equation}\label{MA'}
		M_{E,N}^{-1}\op_N(a)M_{E,N}=\op_N(a \circ E).
	\end{equation}

	We quote from ~\cite[(2.47)]{DJ23} the following adjoint formula
	\begin{equation}\label{3.15f}
		\left(\op_{N}(a)\right)^*=\op_{N}(\overline{a}).
	\end{equation}

	For $a\in C^\infty(\bT^{2n})$, we specify our choice of  Fourier transform,
	\begin{equation}
		\widehat{a}(y,\eta) \coloneqq \int_{\bT^{2n}} e^{-2 \pi\sqrt{-1}\langle (y,\eta), ({x},{\xi})\rangle}a({x},{\xi}) d{x}d{\xi}.
	\end{equation}
	We then have the following Fourier decomposition,
	\begin{equation}
		a(z) = \sum_{\omega\in \Z^{2n}} \widehat{a}(\omega)e^{2\pi\sqrt{-1} \langle z,\omega\rangle}.
	\end{equation}
	Using~\eqref{3.3}, we quantize the Fourier decomposition,
	\begin{equation}\label{3.18f}
		\op_{N}(a) = \sum_{\omega \in \Z^{2n}} \widehat{a}(\omega)U_{\omega/N},
	\end{equation}
	where $\sum_{\omega \in \Z^{2n}} |\widehat{a}(\omega)| < \infty$ by the smoothness of $a$.

	\subsection{Pseudolocality property}\label{s7.3}
	
	We now build on Proposition~\ref{p3.2} to prove our full pseudolocality result. Recall the $B$-block notation from~\eqref{3.5x}.
	
	\begin{prop}\label{prop:pseudolocality}
		Let $E,E'\in\Sp_{2n}(\Z)$ and $a, b \in C^\infty(\bT^{2n})$. Assume that the $B$-block of $EE'$ satisfies $\det B\neq0$. Then
		\begin{equation}\label{3.28}
			\begin{split}
				&\frac{1}{\dim_\C(\cH_N)}\left|\Tr^{\cH_N}[M_{E,N} \op_{N}(a)\op_N(b)  M_{E',N}]\right|\\ &\leqslant \frac{|\det B|^{\frac{1}{2}}}{N^{\frac{n}{2}}} \Blk\sum_{\omega \in \Z^{2n}} |\widehat{a}(\omega)|\Brk\Blk\sum_{\omega \in \Z^{2n}} |\widehat{b}(\omega)|\Brk.
			\end{split}
		\end{equation}
	\end{prop}

	\begin{pro}
		From~\eqref{3.3x},~\eqref{3.8xx}, and~\eqref{3.18f}, we have
		\begin{equation}\label{3.25}
			\begin{split}
				& \frac{1}{\dim_\C(\cH_N)}\left|\Tr^{\cH_N}[M_{E,N} \op_N(a) \op_N(b) M_{E',N}] \right| \\
				&= 
				\frac{1}{\dim_\C(\cH_N)}\bbv\sum_{\omega, \omega' \in \Z^{2n}} \widehat{a}(\omega) \widehat{b}(\omega')\sum_{j \in \Z^n_N} \langle M_{E,N}U_{\omega/N}U_{\omega'/N} M_{E',N} e_j, e_j\rangle \bbv\\
				&\leqslant \frac{1}{\dim_\C(\cH_N)} \Big(\sum_{\omega, \omega' \in \Z^{2n}} |\widehat{a}(\omega)|| \widehat{b}(\omega')|\Big) \\
				&\qquad\qquad\qquad\cdot\Big(\sup_{\omega, \omega' \in \Z^{2n}}\Bv\sum_{j \in \Z^n_N} \langle M_{EE',N}U_{(E')^{T}\frac{\omega+\omega'}{N}}  e_j, e_j\rangle \Bv\Big).
			\end{split}
		\end{equation}
		We have that  $(E')^{T}\frac{\omega+\omega'}{N} = \frac{1}{N}\omega''$ for some $\omega''=(y'',\eta'')\in\Z^{2n}$. By~\eqref{3.19x}, 
		\begin{equation}\label{6.14}
			\begin{split}
			\Bv\sum_{j \in \Z^n_N} \langle M_{EE',N}U_{(E')^{T}\frac{\omega+\omega'}{N}}  e_j, e_j\rangle \Bv& =\bbv\sum_{j \in \Z^n_N} e^{\frac{2 \pi\sqrt{-1}}{N} \langle y'', j\rangle}\langle M_{EE',N}   e_{j-\eta''}, e_j\rangle \bbv\\
				&\leqslant \sum_{j \in \Z^n_N} \bv\langle M_{EE',N}   e_{j-\eta''}, e_j\rangle\bv.
			\end{split}
		\end{equation}
		
		Since the $B$-block of $EE'$ satisfies $\det B\neq0$, by Proposition~\ref{p3.2},~\eqref{3.25}, and~\eqref{6.14}, we conclude~\eqref{3.28}.\qed
	\end{pro}

	\subsection{Trace asymptotics for quantum observables}\label{s7.4}
	
Finally, we have the following trace asymptotics. For $a \in C^\infty(\bT^{2n})$, we use the norm 
	\begin{equation}
		\lV a\rV_{C^k} \coloneqq \sup_{({x}, {\xi})\in\bT^{2n}} \max_{|\alpha| \leqslant k} \lv\partial^\alpha a({x}, {\xi})\rv.
	\end{equation} 
	\begin{prop}
		For all $a,b \in C^\infty(\bT^{2n})$, for any $k \geqslant 2n+1$,
		\begin{equation}\label{7.16x}
			\begin{split}
				&\frac{1}{\dim_\C(\cH_N)}\Tr^{\cH_N}[\op_N(a) \op_N(b)]\\
				&=\int_{\bT^{2n}} a(x,\xi) b(x,\xi)dxd\xi+O_{n,k}\Big(N^{-k} \sum_{\substack{0 \leqslant k_1, k_2 \leqslant k \\ k_1 + k_2 =k}} \lV a\rV_{C^{k_1}}\lV b\rV_{C^{k_2}}\Big).
			\end{split}
		\end{equation}
	\end{prop}
	
	\begin{pro}
		1. Using~\eqref{3.18f}, we know that
		\begin{equation}\label{3.35}
			\begin{split}
				&\frac{1}{\dim_\C(\cH_N)}\Tr^{\cH_N}[\op_N(a) \op_N(b)]\\
				&=\frac{1}{\dim_\C(\cH_N)} \sum_{\omega, \omega' \in \Z^{2n}} \widehat{a}(\omega) \widehat{b}(\omega') \Tr^{\cH_N}[U_{\omega/N} U_{\omega'/N}].
			\end{split}
		\end{equation}
		2. We claim that for $\omega=(y,\eta), \omega'=(y',\eta')\in\Z^{2n}$,
		\begin{equation}\label{3.44f}
			\begin{split}
				&\frac{1}{\dim_\C(\cH_N)}\Tr^{\cH_N}[U_{\omega/N}U_{\omega'/N} ]\\
				&= \begin{cases}
					(-1)^{\sigma(\omega,\omega'') + N \langle y'',\eta''\rangle}, & \omega' = -\omega+ N\omega'',\\
					0, & \omega' \not\equiv -\omega\bmod N.
				\end{cases}
			\end{split}
		\end{equation}
		
		Since
		\begin{equation}\label{eq:first-trace}
			\frac{1}{\dim_\C(\cH_N)}\Tr^{\cH_N}[U_{\omega/N}U_{\omega'/N} ] = \frac{1}{N^n}\sum_{j \in \Z^n_N} \langle U_{\omega/N} U_{\omega'/N} e_j,e_j\rangle,
		\end{equation}
		the $\omega' = -\omega+ N\omega''$ case  follows directly from~\eqref{3.3x} and~\eqref{3.19x}.
		
		We now consider the case where $\omega  \not\equiv -\omega' \bmod N$. By~\eqref{3.19x} and~\eqref{eq:first-trace}, we have
		\begin{equation}
			\begin{split}
				&\frac{1}{\dim_\C(\cH_N)}\Tr^{\cH_N}[U_{\omega/N}U_{\omega'/N} ]\\
				&= \frac{e^{-\frac{\pi\sqrt{-1}}{N} (\langle y,\eta\rangle + \langle y', \eta' \rangle)}}{N^n} \sum_{j \in \Z^n_N} e^{\frac{2 \pi\sqrt{-1} }{N} \langle y+y', j \rangle} \langle e_{j-\eta' - \eta}, e_{j}\rangle.
			\end{split}
		\end{equation}
		Clearly, this is zero if $-\eta \not\equiv \eta' \bmod N$. So we assume that $-\eta\equiv\eta'\bmod N$ and $-y \not\equiv y' \bmod N$. 
		Then,
		\begin{equation}\label{3.38}
			\begin{split}
				&\frac{1}{\dim_\C(\cH_N)}\Tr^{\cH_N}[U_{\omega/N}U_{\omega'/N} ]\\
				&= \frac{e^{-\frac{\pi\sqrt{-1}}{N} (\langle y, \eta\rangle + \langle y',\eta' \rangle)}}{N^n}\sum_{j \in \Z^n_N}  e^{\frac{2 \pi\sqrt{-1} }{N} \langle y+y', j \rangle} = 0,
			\end{split}
		\end{equation}
		which completes the proof of~\eqref{3.44f}.
		
		3. By~\eqref{3.35} and~\eqref{3.44f}, we see that
		\begin{equation}\label{3.47}
			\begin{split}
				&\frac{1}{\dim_\C(\cH_N)}\Tr^{\cH_N}[\op_N(a) \op_N(b)]\\
				&=\sum_{\omega'\in\Z^{2n}}\sum_{\omega \in \Z^{2n}} (-1)^{\sigma(\omega,\omega')+c(\omega')}\widehat{a}(-\omega)\widehat{b}(\omega-N\omega'),
			\end{split}
		\end{equation}
		where $c(\omega')$ is an integer that depends only on $\omega'$ with $c(0)=0$.

		By the smoothness of $a,b$, we know $\sum_{\omega \in \Z^{2n}}\vert\widehat{a}(\omega)\vert, \sum_{\omega \in \Z^{2n}}\vert\widehat{b}(\omega)\vert< \infty$.
		Therefore, by Parseval's identity, the $\omega'=0$ term in~\eqref{3.47} is
		\begin{equation}\label{3.48f}
			\sum_{\omega \in \Z^{2n}} \widehat{a}(-\omega) \widehat{b}(\omega)=\int_{\bT^{2n}} a({x},{\xi}) b({x},{\xi}) d{x}d{\xi}.
		\end{equation}

		4. We finish our estimate of~\eqref{3.47} by showing that the sum over $\omega'\in \Z^{2n}\setminus \{0\}$ gives the error term.
		We have
		\begin{equation}\label{3.49f}
			\begin{split}
				&\sum_{\omega' \in \Z^{2n} \setminus \{0\}} \sum_{\omega \in \Z^{2n}} (-1)^{\sigma(\omega, \omega')+c(\omega')}\widehat{a}(-\omega) \widehat{b}(\omega-N\omega')\\
				&=\sum_{\omega'\in\Z^{2n} \setminus \{0\}} \sum_{\omega, \omega'' \in \Z^{2n}} (-1)^{\sigma(\omega,\omega') + c(\omega')} \widehat{a}(-\omega) \widehat{b}(\omega'')  \int_{\bT^{2n}} e^{2 \pi\sqrt{-1} \langle \omega''-\omega +N\omega',z\rangle}dz \\
				&= \sum_{\omega'\in\Z^{2n}\setminus \{0\}}  \int_{\bT^{2n}}  \Big(\sum_{\omega\in \Z^{2n}} (-1)^{\sigma(\omega,\omega') + c(\omega')}  \widehat{a}(-\omega)e^{2 \pi\sqrt{-1} \langle-\omega,z\rangle}\Big)\\
				&\qquad\qquad\qquad\qquad\cdot\Big(\sum_{\omega''\in \Z^{2n}} \widehat{b}(\omega'')e^{2 \pi\sqrt{-1} \langle \omega'' +N\omega',z\rangle}\Big)dz\\
				&=\sum_{\omega'\in \Z^{2n}\setminus \{0\}} (-1)^{c(\omega')}  \int_{\bT^{2n}}  \Big(\sum_{\omega\in \Z^{2n}} (-1)^{\sigma(\omega,\omega') }  \widehat{a}(-\omega)e^{2 \pi\sqrt{-1} \langle-\omega,z\rangle}\Big)\\
				&\qquad\qquad\qquad\qquad\qquad\qquad\quad\cdot b(z)e^{2 \pi N\sqrt{-1}\langle \omega',z\rangle}dz.
			\end{split}
		\end{equation}
		
		We note that
		\begin{equation}
			\sum_{\omega\in \Z^{2n}} (-1)^{\sigma(\omega,\omega') }  \widehat{a}(-\omega)e^{2 \pi\sqrt{-1}\langle-\omega,z\rangle}= a\left(({x},{\xi})-\tfrac{(\eta',y')}{2}\right).
		\end{equation}
		
		We then use an integration by parts argument to see that for $\omega'\not\equiv 0$ and all $k \geqslant 0$,
		\begin{equation}
			\begin{split}
				\left\vert\int_{\bT^{2n}}a(({x}, {\xi})-(\eta',y')/2)b(z) e^{2 \pi\sqrt{-1} \langle N\omega',z\rangle}dxd\xi\right\vert \\
				\leqslant C_k (2 \pi N)^{-k} \lV \omega'\rV_\infty^{-k} \lV a\rV_{C^k}\lV b\rV_{C^k}.
			\end{split}
		\end{equation}
		Therefore,
		\begin{equation}\label{6.26}
			\begin{split}
				&\sum_{\omega'\in \Z^{2n} \setminus \{0\}}\lv\int_{\bT^{2n}}a(({x}, {\xi})-(\eta',y')/2)b(z)e^{2 \pi\sqrt{-1} \langle N\omega',z\rangle}dxd\xi\rv\\
				&\leqslant C_{n,k} N^{-k} \sum_{\substack{0 \leqslant k_1, k_2 \leqslant k,\\ k_1 + k_2 =k}} \lV a\rV_{C^{k_1}}\lV b\rV_{C^{k_2}}  \sum_{\omega'\in\Z^{2n} \setminus \{0\}}\lV \omega'\rV_\infty^{-k}\\
				&\leqslant C_{n,k} N^{-k} \sum_{\substack{0 \leqslant k_1, k_2 \leqslant k,\\ k_1+k_2=k}} \lV a\rV_{C^{k_1}}\lV b\rV_{C^{k_2}} \sum_{\ell =1}^\infty \sum_{\lV \omega'\rV_{\infty} =\ell} \ell^{-k}.
			\end{split}
		\end{equation}		
		Since
		\begin{equation}
			\lv\{\omega' \in \Z^{2n}\mid \lV \omega'\rV_\infty = \ell\}\rv \leqslant C_n \ell^{2n-1},
		\end{equation}
		the final sum in~\eqref{6.26} is bounded if $k \geqslant(2n+1)$, which completes the proof.\qed
	\end{pro}

	Setting $b=1$, we get the following result. In~\eqref{7.16x}, the $k$-dependence for the implied constant came from taking derivatives of the product $ab$ in the integration by parts step.  Thus, the implied constant in the following no longer has any $k$-dependence. 
	\begin{prop}
		For all $a \in C^\infty(\bT^{2n})$, for any $k \geqslant (2n+1)$,
		\begin{equation}
			\frac{1}{\dim_\C(\cH_N)}\Tr^{\cH_N}[\op_N(a)] = \int_{\bT^{2n}}a(z)dz+O_n\big(2\pi N)^{-k}\lV a\rV_{C^{k}}\big),
		\end{equation}
		where the constant in $O(\cdot)$ depends only on $n$.
	\end{prop}

	\section{Mixed Quantization}\label{S8}

	In this section, we combine the kernel operators introduced in \S\,\ref{S6} with the quantization of functions on the torus introduced in \S\,\ref{S7} to construct our discrete mixed quantization. In \S\,\ref{s9.1}, we introduce kernel functions. In \S\,\ref{s8.2}, we study the associated commutator operator. In \S\,\ref{s9.2}, we define the fiberwise quantization of kernel functions as discrete mixed quantization.
	
	For mixed quantization on manifolds, we refer to ~\cite[\S\,4]{MR4808253} and ~\cite[\S\,3]{ovadia2025mixedquantizationpartialhyperbolicity}. For its discrete counterpart, we refer to ~\cite[\S\,10]{disc}.

	\subsection{Kernel functions}\label{s9.1}

	Similarly to~\eqref{7.58} and~\eqref{7.59}, we define
	\begin{equation}\label{9.1}
		\begin{split}
			\mathscr{X} & \coloneqq \Gamma\backslash\big(\mathbf{T}_d\times \mathbb{T}^{2n}\big)\\
			&=\big\{(\bm{v},z)\in \mathbf{T}_d\times \mathbb{T}^{2n}\big\}\big/(\bm{v},z)\sim(\gamma \bm{v},\chi(\gamma) z) \text{ for every }\gamma\in\Gamma.
		\end{split}
	\end{equation}
	Then the space of functions and the Laplacian are given by
	\begin{equation}\label{9.2}
		\begin{split}
			C^\infty(\mathscr{X})& \coloneqq C^\infty_\Gamma(\mathbf{T}_d\times \mathbb{T}^{2n})\\
			&=\big\{\mathscr{U}\mid \mathscr{U}(\bm{v},z)=\mathscr{U}(\gamma \bm{v},\chi(\gamma) z)\text{ for every }(\bm{v},z)\in \mathbf{T}_d\times \mathbb{T}^{2n}, \gamma\in\Gamma\big\},\\
			(\Delta^{\mathscr{X}}\mathscr{U})(\bm{v},z)& \coloneqq \sum_{\bm{v}'\sim \bm{v}}\mathscr{U}(\bm{v}',z).
		\end{split}
	\end{equation}

	Similarly to~\eqref{7.1n} and~\eqref{7.2n}, we define
	\begin{equation}
		\begin{split}
			\mathscr{K}& \coloneqq \Gamma\backslash\big(\mathbf{T}_d^{\times 2}\times \mathbb{T}^{2n}\big)\\
			&= \big\{(\bm{v},\bm{v}',z)\in \mathbf{T}_d^{\times2}\times \mathbb{T}^{2n}\big\}\big/(\bm{v},\bm{v}',z)\sim(\gamma \bm{v},\gamma \bm{v}',\chi(\gamma)z) \text{ for every }\gamma\in\Gamma.
		\end{split}
	\end{equation}
	Then the space of smooth functions is given by
	\begin{equation}\label{xx8.4}
		\begin{split}
			C^\infty(\mathscr{K}) \coloneqq \big\{\mathscr{Q}\colon\mathbf{T}_d^{\times 2}\times \mathbb{T}^{2n}\to \mathbb{C}\mid\ & \mathscr{Q}(\bm{v},\bm{v}',z)=\mathscr{Q}(\gamma \bm{v},\gamma \bm{v}',\chi(\gamma) z)\\
			&\text{for every }(\bm{v},\bm{v}',z)\in \mathbf{T}_d^{\times 2}\times \mathbb{T}^{2n}, \gamma\in\Gamma\big\},
		\end{split}
	\end{equation}
	and we call its elements \emph{kernel functions}. We define $\mathscr{K}_{(k)}, \mathscr{K}_{(\leqslant k)}$ analogously to $K_{(k)}$, $K_{(\leqslant k)}$ in ~\eqref{7.4}, and $\mathscr{Q}_{(k)}$ analogously to $Q_{(k)}$ in~\eqref{7.5nnn}.

	Similarly to~\eqref{7.7n}, we define the normalized $L^{2}$-norm on $C^\infty(\mathscr{K})$ by
	\begin{equation}\label{9.7}
		\lV\mathscr{Q}\rV_{L^2(\mathscr{K})}^2=\frac{1}{\lv D\rv}\sum_{\bm{v}\in D,\bm{v}'\in \mathbf{T}_d}\lV \mathscr{Q}(\bm{v},\bm{v}',\cdot)\rV_{L^2(\mathbb{T}^{2n})}^2.
	\end{equation}

	\subsection{Commutator operators}\label{s8.2}

	As in~\eqref{7.27n}, we define
	\begin{equation}\label{9.9}
		\begin{split}
			&\mathrm{ad}_{\Delta^{\mathscr{X}}}\colon C^{\infty}(\mathscr{K}_{(k)}) \to C^{\infty}(\mathscr{K}_{(k-1)}\cup \mathscr{K}_{(k+1)}),\\
			&\mathrm{ad}_{\Delta^{\mathscr{X}}}(\mathscr{Q})(\bm{v},\bm{v}',z)=\sum_{\bm{v}''\sim \bm{v} }\mathscr{Q}(\bm{v}'',\bm{v}',z)-\sum_{\bm{v}'''\sim \bm{v}'}\mathscr{Q}(\bm{v},\bm{v}''',z).
		\end{split}
	\end{equation}

	Similarly to~\eqref{7.5nn}, we define $\mathbbm{1}_{(k)}\in C^\infty(\mathscr{K}_{(k)})$ by
	\begin{equation}\label{9.5}
		\mathbbm{1}_{(k)}(\bm{v},\bm{v}',z) \coloneqq \mathbbm{1}_{\{d_{\mathbf{T}_d}(\bm{v},\bm{v}')=k\}}(\bm{v},\bm{v}'),
	\end{equation}
	which depends only on the distance in the base $K$.
	
We next introduce some dynamical results, for which it is more natural to work in $L^2(\mathscr{K})$ rather than in $C^\infty(\mathscr{K})$. This is analogous to the proof of the classical quantum ergodicity, treating smooth observables as square integrable functions when applying dynamical results.

	\begin{prop}
		Assume  \textbf{[ERG]}~\eqref{5.11}. Then on $L^2(\mathscr{K})$,
		\begin{equation}\label{9.10}
			\ker(\mathrm{ad}_{\Delta^{\mathscr{X}}})=\mathrm{span}\{\mathbbm{1}_{(k)}\mid k\in\mathbb{N}\},
		\end{equation}
		where the span is understood as the closed linear span.
	\end{prop}
	
	\begin{pro}
		The proof of Proposition~\ref{L7.17} generalizes to the current setting without essential modifications. Therefore,
		\begin{equation}
			\begin{split}
				&\ker(\mathrm{ad}_{\Delta^{\mathscr{X}}})\\
				&=\mathrm{span}\big\{f\cdot \mathbbm{1}_{(k)}\mid k\in\mathbb{N},f\in L^2(\mathbb{T}^{2n}),f\circ \chi(\gamma)=f\text{ for every }\gamma\in\Gamma\big\}.
			\end{split}
		\end{equation}
		By \textbf{[ERG]}~\eqref{5.11}, the only $\chi(\Gamma)$-invariant functions on $\mathbb{T}^{2n}$ are constant functions.\qed
	\end{pro}

	We define the normalized average of $\mathscr{Q}\in L^2(\mathscr{K}_{(k)})$ by
	\begin{equation}\label{9.12}
		\begin{split}
			\langle \mathscr{Q}\rangle_{\mathscr{K}}& \coloneqq \frac{1}{\lv \pa B_{\mathbf{T}_d}(o,k)\rv}\big\langle \mathscr{Q},\mathbbm{1}_{(k)}\big\rangle_{L^2(\mathscr{K}_{(k)})}\\
			&=\frac{1}{\lv D\rv\cdot \lv \pa B_{\mathbf{T}_d}(o,k)\rv}\sum_{\bm{v}\in D,\bm{v}'\in \mathbf{T}_d}\int_{\mathbb{T}^{2n}}\mathscr{Q}(\bm{v},\bm{v}',z)dz.
		\end{split}
	\end{equation}

	\begin{prop}
		Assume  \textbf{[ERG]}~\eqref{5.11}. Then the orthogonal projection
		\begin{equation}\label{8.10}
			P^{\ker(\mathrm{ad}_{\Delta^{\mathscr{X}}})}\colon L^2(\mathscr{K})\to \ker(\mathrm{ad}_{\Delta^{\mathscr{X}}})
		\end{equation}
		can be expressed by
		\begin{equation}\label{9.13}
			P^{\ker(\mathrm{ad}_{\Delta^{\mathscr{X}}})}\mathscr{Q}=\langle \mathscr{Q}\rangle_{\mathscr{K}}\mathbbm{1}_{(k)}.
		\end{equation}
	\end{prop}
	
	\begin{pro}
		This follows directly from~\eqref{9.7} and~\eqref{9.10}.\qed
	\end{pro}

	\begin{prop}
	Assume \textbf{[ERG]}~\eqref{5.11}. Let $\mathscr{Q}\in L^2(\mathscr{K})$ satisfy $\mathscr{Q}\perp \mathrm{span}\{\mathbbm{1}_{(k)}\mid k\in\mathbb{N}\}$. Then,
		\begin{equation}\label{9.16}
			\lim_{t_0\to\infty}\bbV\frac{1}{{t_0}}\int_0^{t_0}e^{\sqrt{-1}t\mathrm{ad}_{\Delta^{\mathscr{X}}}}\mathscr{Q}dt\bbV_{L^2(\mathscr{K})}=0.
		\end{equation}
	\end{prop}
	
	\begin{pro}
		Analogously to~\eqref{7.36}, $\mathrm{ad}_{\Delta^{\mathscr{X}}}$ is selfadjoint. The proof then follows directly from~\eqref{9.10} and the von Neumann mean ergodic theorem.\qed
	\end{pro}

	\subsection{The fiberwise quantization of kernel functions}\label{s9.2}

	Recall from~\eqref{4.2x} and ~\eqref{4.5ns} that $F_N$ is induced by the unitary representation $\rho_N$ of $\Gamma$ on $\mathcal{H}_N$, and from ~\eqref{7.2n} that the space of kernel operators is given by
	\begin{equation}\label{8.14}
		\begin{split}
			C(K, K^{\mathrm{End}(F_N)})=\big\{Q\colon\mathbf{T}_d^{\times 2}\to \mathrm{End}(\mathcal{H}_N) \mid\ &Q(\gamma \bm{v},\gamma  \bm{v}')=\rho_N(\gamma)Q(\bm{v},\bm{v}')\rho_N(\gamma)^{-1}\\
			&\text{for every } (\bm{v},\bm{v}')\in \mathbf{T}_d^{\times 2}, \gamma\in\Gamma\big\}.
		\end{split}
	\end{equation}
	Since $\rho_N(\gamma)$ differs from $M_{\chi(\gamma),N}$ only by a unit complex factor, the same equivariance identity holds with $\rho_N(\gamma)$ replaced by $M_{\chi(\gamma),N}$. We now construct such kernel operators by quantizing kernel functions in the torus variable.

	Consider a kernel function $\mathscr{Q}\in C^\infty(\mathscr{K})$, defined in~\eqref{xx8.4}. For any $(\bm{v},\bm{v}')\in \mathbf{T}_d^{\times 2}$, we have $\mathscr{Q}(\bm{v},\bm{v}',\cdot)\in C^\infty(\mathbb{T}^{2n})$. Applying the quantization~\eqref{7.7} to this function, we get $\op_{N}\big(\mathscr{Q}(\bm{v},\bm{v}',\cdot)\big)\in\mathrm{End}(\mathcal{H}_N)$. This defines an endomorphism-valued function $\op_{N}(\mathscr{Q})$ on $\mathbf{T}_d^{\times 2}$ by
	\begin{equation}
		\op_{N}(\mathscr{Q})(\bm{v},\bm{v}')=\op_{N}\big(\mathscr{Q}(\bm{v},\bm{v}',\cdot)\big).
	\end{equation}
	From~\eqref{MA'} and~\eqref{xx8.4}, we obtain
	\begin{equation}\label{10.30}
		\begin{split}
			\op_{N}(\mathscr{Q})(\gamma \bm{v},\gamma \bm{v}')&=\op_{N}\big(\mathscr{Q}(\gamma \bm{v},\gamma \bm{v}',\cdot)\big)\\
			&=\op_{N}\big(\mathscr{Q}(\bm{v},\bm{v}',\chi(\gamma^{-1})\cdot)\big)\\
			&=M_{\chi(\gamma),N}\op_{N}\big(\mathscr{Q}(\bm{v},\bm{v}',\cdot)\big)M_{\chi(\gamma),N}^{-1}\\
			&=M_{\chi(\gamma),N}\op_{N}(\mathscr{Q})(\bm{v},\bm{v}')M_{\chi(\gamma),N}^{-1}.
		\end{split}
	\end{equation}
	Comparing~\eqref{8.14} with~\eqref{10.30}, we see that $\op_{N}(\mathscr{Q})$ is a kernel operator. 
	
	Since the quantization acts only in the torus variable, it preserves the propagation of kernel functions. In summary, we obtain the fiberwise geometric quantization map
	\begin{equation}\label{9.17}
		\begin{split}
			&\op_{N}\colon C^\infty(\mathscr{K}_{(k)})\to C\big(K_{(k)},K^{\mathrm{End}(F_N)}\big),\\
			&\op_{N}(\mathscr{Q})(\bm{v},\bm{v}')=\op_{N}(\mathscr{Q}(\bm{v},\bm{v}',\cdot))
		\end{split}
	\end{equation}
	Kernel operators provide the semiclassical structure in the graph variable, while quantization acts on the torus variable. The fiberwise quantization~\eqref{9.17} combines these two structures, giving a natural model of \emph{discrete mixed quantization} and motivating the title of this paper.
	
	The following result expresses the compatibility of fiberwise quantization with the commutator operators.
	\begin{prop}
		Let $\mathscr{Q}\in C^\infty(\mathscr{K})$. Then for $N\in 2\mathbb{N}$,
		\begin{equation}\label{9.18}
			\mathrm{ad}_{\Delta^{F_N}}\op_{N}(\mathscr{Q})=\op_{N}(\mathrm{ad}_{\Delta^{\mathscr{X}}}\mathscr{Q}).
		\end{equation}
	\end{prop}
	
	\begin{pro}
		This is immediate from~\eqref{7.27n},~\eqref{9.9}, and~\eqref{9.17}.\qed
	\end{pro}

	\begin{prop}\label{p9.7}
		Assume \textbf{[BIV]}~\eqref{5.10}. Let $\mathscr{Q}\in C^\infty(\mathscr{K}_{(\leqslant k)})$. Then as $N\in2\mathbb{N}$ tends to infinity,
		\begin{equation}\label{9.20}
			\lV\op_{N}(\mathscr{Q})\rV_{\mathrm{HS}(X,F_N)}^2=\lV\op_{N}(\mathscr{Q})\rV_{L^2(K_{(\leqslant k)},K^{\mathrm{End}(F_N)})}^2+o(1).
		\end{equation}
	\end{prop}
	
	\begin{pro}
		By~\eqref{7.10n}, it suffices to show that
		\begin{equation}\label{9.21}
			\begin{split}
				\lim_{N\to\infty}\sum_{\substack{\bm{v},\bm{v}'\in D,\  \gamma,\gamma'\in\Gamma\\ \gamma\neq \gamma'}}\frac{\mathrm{Tr}^{\mathcal{H}_N}}{\dim_{\mathbb{C}}\mathcal{H}_N}\Big[M_{\chi((\gamma')^{-1}),N}\op_{N}(\overline{\mathscr{Q}})(\bm{v},\gamma'\bm{v}')&\\
				\op_{N}(\mathscr{Q})(\bm{v},\gamma \bm{v}')M_{\chi(\gamma),N}&\Big]=0.
			\end{split}
		\end{equation}
		We note that the sum in~\eqref{9.21} is finite and  $(\gamma')^{-1}\gamma$ is nonidentity for $\gamma\neq\gamma'$. By \textbf{[BIV]}~\eqref{5.10} and the pseudolocality~\eqref{3.28}, each summand tends to zero as $N\to\infty$. This proves~\eqref{9.21}, and hence~\eqref{9.20}.\qed
	\end{pro}

	\begin{prop}\label{l9.6n}
		Let $\mathscr{Q}\in C^\infty(\mathscr{K}_{(\leqslant k)})$. Then as $N\in2\mathbb{N}$ tends to infinity,
		\begin{equation}\label{9.19}
			\lV\op_{N}(\mathscr{Q})\rV_{L^2(K,K^{\mathrm{End}(F_N)})}^2=\lV \mathscr{Q}\rV_{L^2(\mathscr{K})}^2+o(1).
		\end{equation}
	\end{prop}

	\begin{pro}
		This is a direct consequence of~\eqref{7.7n},~\eqref{7.16x}, and~\eqref{9.7}.\qed
	\end{pro}

	\section{Quantum Ergodicity}\label{s6n}
	
	In this section, we prove our quantum ergodicity result. In \S\,\ref{s10.1}, we define the quantum variance and establish its basic properties. In \S\,\ref{s9.2x}, we combine the pseudolocality, ergodicity, and mixed quantization results to prove eigensection equidistribution.

For related results on discrete quantum ergodicity, we refer to ~\cite[Theorem 1.1]{AlM15},~\cite[Theorem 1]{MR3567266}, and ~\cite[Theorem 11.3]{disc}.

	\subsection{Quantum variance}\label{s10.1}
	
	We define the normalized quantum variance of a kernel operator $Q\in C(K,K^{\mathrm{End}(F_N)})$ from~\eqref{8.14} by
	\begin{equation}\label{9.1x}
		\mathrm{Var}(Q) \coloneqq \frac{1}{\dim_{\mathbb{C}}F_N}\sum_{i=1}^{\dim_{\mathbb{C}}F_N}\bv\big\langle Qu_{N,i},u_{N,i}\big\rangle_{L^2(X,F_N)}\bv^2.
	\end{equation}
	
	Since $\Delta^{F_N}$ is selfadjoint and $(u_{N,i})$ consists of eigensections, we have for any $j\in\mathbb{N}$ and $t_0>0$,
	\begin{equation}\label{10.4}
		\mathrm{Var}(Q) =\mathrm{Var}\bigg(\sum_{i\leqslant j}\frac{1}{{t_0}}\int_0^{t_0}\frac{(t\sqrt{-1}\mathrm{ad}_{\Delta^{F_N}})^i}{i!}Qdt\bigg).
	\end{equation}
	Indeed, the operator inside the variance on the right can be written as $(Q+\mathrm{ad}_{\Delta^{F_N}}Q')$ for some $Q'\in C(K,K^{\mathrm{End}(F_N)})$. For each eigensection $u_{N,i}$,
	\begin{equation}
		\begin{split}
			&\big\langle \big(\mathrm{ad}_{\Delta^{F_N}}Q'\big)u_{N,i},u_{N,i}\big\rangle_{L^2(X,F_N)}\\
			&=\big\langle Q'u_{N,i},\Delta^{F_N}u_{N,i}\big\rangle_{L^2(X,F_N)}-\big\langle Q'\Delta^{F_N}u_{N,i},u_{N,i}\big\rangle_{L^2(X,F_N)}=0.
		\end{split}
	\end{equation}
	This gives~\eqref{10.4}.

	Since $(u_{N,i})$ is an orthonormal basis, we also have
	\begin{equation}\label{10.2}
		\mathrm{Var}(Q)\leqslant\frac{1}{\dim_{\mathbb{C}}F_N}\sum_{i=1}^{\dim_{\mathbb{C}}F_N}\lV Qu_{N,i}\rV_{L^2(X,F_N)}^2=\lV Q\rV_{\mathrm{HS}(X,F_N)}^2.
	\end{equation}

	\subsection{Eigensection equidistribution}\label{s9.2x}

	We now state and prove our quantum ergodicity result. In particular, when $k=0$ in the following theorem, it is equivalent to  Theorem~\ref{C9'} by Proposition~\ref{lem:unitary-equiv}. Recall the polynomial $h_k$ defined in~\eqref{6.31}, the normalized average $\langle\mathscr{Q}_{(k)}\rangle_{\mathscr{K}}$ defined in~\eqref{9.12}, and the fiberwise quantization $\op_N(\mathscr{Q}_{(k)})$ defined in~\eqref{9.17}.

	\begin{theo}\label{t9.1}
		In the setting of~\eqref{n2.6}, assume \textbf{[BIV]}~\eqref{5.10} and \textbf{[ERG]}~\eqref{5.11}. Let $\mathscr{Q}_{(k)}\in C^\infty(\mathscr{K}_{(k)})$. Then
		\begin{equation}\label{10.11}
			\begin{split}
				\lim_{\substack{N\in2\mathbb{N},\\ N\to\infty}}\frac{1}{\dim_{\mathbb{C}}F_N}\sum_{i=1}^{\dim_{\mathbb{C}}F_N}\bv\big\langle \op_{N}(\mathscr{Q}_{(k)})u_{N,i},u_{N,i}\big\rangle_{L^2(X,F_N)}&\\
				-\langle \mathscr{Q}_{(k)}\rangle_{\mathscr{K}} h_k(\lambda_{N,i})&\bv^2=0.
			\end{split}
		\end{equation}
	\end{theo}

	\begin{pro}
		By~\eqref{7.5nn},~\eqref{6.31},~\eqref{9.5},~\eqref{9.17}, and the fact that the quantization of the constant function $\mathbbm{1}_{\mathbb{T}^{2n}}$ on $\mathbb{T}^{2n}$ is the identity on $\mathcal{H}_N$, we get
		\begin{equation}
			\op_N(\mathbbm{1}_{(k)})=\mathrm{Id}_{(k)}=h_k(\Delta^{F_N}).
		\end{equation}
		Therefore,
		\begin{equation}
			\big\langle\op_N(\mathbbm{1}_{(k)})u_{N,i},u_{N,i}\big\rangle_{L^2(X,F_N)}=h_k(\lambda_{N,i}),
		\end{equation}
		and by~\eqref{9.1x}, the limit in~\eqref{10.11} is
		equivalent to showing
		\begin{equation}
			\lim_{\substack{N\in2\mathbb{N},\\ N\to\infty}}\mathrm{Var}\big(
			\op_N\big(\mathscr{Q}_{(k)}-\langle\mathscr{Q}_{(k)}\rangle_{\mathscr{K}}\mathbbm{1}_{(k)}\big)\big)=0.
		\end{equation}
		
		Without loss of generality, we may suppose that $\mathscr{Q}_{(k)}\perp \mathbbm{1}_{(k)}$. From ~\eqref{10.4} and~\eqref{10.2}, we get
		\begin{equation}\label{10.12}
			\begin{split}
				&\lim_{\substack{N\in2\mathbb{N},\\ N\to\infty}}	\mathrm{Var}\big(\op_{N}(\mathscr{Q}_{(k)})\big)\\
				&=\lim_{\substack{N\in2\mathbb{N},\\ N\to\infty}}	\mathrm{Var}\Big(\sum_{i\leqslant j}\frac{1}{{t_0}}\int_0^{t_0}\frac{(t\sqrt{-1}\mathrm{ad}_{\Delta^{F_N}})^i}{i!}\op_{N}(\mathscr{Q}_{(k)})dt\Big)\\
				&\leqslant\lim_{\substack{N\in2\mathbb{N},\\ N\to\infty}} \bbV\sum_{i\leqslant j}\frac{1}{{t_0}}\int_0^{t_0}\frac{(t\sqrt{-1}\mathrm{ad}_{\Delta^{F_N}})^i}{i!}\op_{N}(\mathscr{Q}_{(k)})dt\bbV_{\mathrm{HS}(X,F_N)}^2.
			\end{split}
		\end{equation}
		By~\eqref{9.18}, the limit~\eqref{10.12} is equal to
		\begin{equation}\label{9.7x}
			\lim_{\substack{N\in2\mathbb{N},\\ N\to\infty}}\BV \op_{N}\Big(\sum_{i\leqslant j}\frac{1}{{t_0}}\int_0^{t_0}\frac{(t\sqrt{-1}\mathrm{ad}_{\Delta^{\mathscr{X}}})^i}{i!}\mathscr{Q}_{(k)}dt\Big)\BV_{\mathrm{HS}(X,F_N)}^2.
		\end{equation}
		By~\eqref{9.20}, under \textbf{[BIV]}~\eqref{5.10}, the limit~\eqref{9.7x} is bounded by
		\begin{equation}\label{10.12n}
			\lim_{\substack{N\in2\mathbb{N},\\ N\to\infty}} \BV \op_{N}\Big(\sum_{i\leqslant j}\frac{1}{{t_0}}\int_0^{t_0}\frac{(t\sqrt{-1}\mathrm{ad}_{\Delta^{\mathscr{X}}})^i}{i!}\mathscr{Q}_{(k)}dt\Big)\BV_{L^2(K,K^{\mathrm{End}(F_N)})}^2.
		\end{equation}
		By~\eqref{9.19}, the limit~\eqref{10.12n} is equal to
		\begin{equation}\label{9.9x}
			\bbV\sum_{i\leqslant j}\frac{1}{{t_0}}\int_0^{t_0}\frac{(t\sqrt{-1}\mathrm{ad}_{\Delta^{\mathscr{X}}})^i}{i!}\mathscr{Q}_{(k)}dt\bbV_{L^2(\mathscr{K})}^2.
		\end{equation}
		Finally, by~\eqref{9.16}, under \textbf{[ERG]}~\eqref{5.11}, we can take $j\to\infty$ and then $t_0\to \infty$ in~\eqref{9.9x} to get~\eqref{10.11}.\qed
	\end{pro}
	
	Note that we cannot simply replace the finite sum in~\eqref{10.12} with the infinite sum $e^{t\sqrt{-1}\mathrm{ad}_{\Delta^{F_N}}}\op_{N}(\mathscr{Q}_{(k)})$, because this would produce a kernel function $e^{t\sqrt{-1}\mathrm{ad}_{\Delta^{\mathscr{X}}}}\mathscr{Q}$ with infinite propagation, and then~\eqref{9.20} would fail.

\section{Applications}\label{examples}

In this section, we give some examples to which Theorems \ref{t5.1} and \ref{t9.1} apply. It suffices to give a representation \eqref{5.7x} satisfying \textbf{[BIV]} \eqref{5.10} and \textbf{[ERG]} \eqref{5.11}.  The construction relies on the classical ping-pong argument, and we mainly refer to de la Harpe~\cite[\S\,\uppercase\expandafter{\romannumeral2}.B]{MR1786869}. We also briefly discuss how the two assumptions can fail.
	
To begin with, for $i\geqslant1$, we consider the matrix
\begin{equation}
E_i=\begin{pmatrix}1+i&i\\
1&1\end{pmatrix}\in\mathrm{SL}_2(\mathbb Z),
\end{equation}
with eigenvalues and eigenvectors
\begin{equation}
\lambda_i^\pm
=\frac{(i+2)\pm\sqrt{i^2+4i}}2,
\qquad x^\pm_i=\begin{pmatrix}\lambda_i^\pm-1\\1\end{pmatrix}.
\end{equation}

Let $E_i$ act on the projective line $\mathbb{RP}^1$. Since $0<\lambda_i^-<1<\lambda_i^+$, we see that $E_i$ expands the line $\mathbb Rx_i^+$ and contracts the line $\mathbb Rx_i^-$. In other words,
\begin{equation}\label{10.3}
	\begin{split}
&\lim_{j\to\infty}E_i^jx=x_i^+,\quad \text{for every }x\neq \mathbb{R}x_i^-,\\
		&\lim_{j\to\infty}E_i^{-j}x=x_i^-, \quad\text{for every }x\neq \mathbb{R}x_i^+.
	\end{split}
\end{equation}

Suppose that $\Gamma$ has free generators $(\gamma_i)_{i=1}^\ell$. Since $\{x^\pm_i\mid 1\leqslant i\leqslant\ell\}$ are pairwise distinct, by \eqref{10.3}, we can choose a sufficiently large $j$ such that the elements $E_i^j$ admit pairwise disjoint neighbor $D_{E_i^{\pm j}}\subset\mathbb{RP}^1$ (neighborhood(s) of $\mathbb{R}x_i^\pm$), all avoiding $\mathbb{R}\big(\begin{smallmatrix}0\\1\end{smallmatrix}\big)$, satisfying the ping-pong relations
	\begin{equation}\label{10.4x}
E\blk\mathbb{RP}^1\setminus D_{E^{-1}}\brk
\subset D_E,\quad\text{for every }E\in\{E_1^{\pm j},\cdots,E_\ell^{\pm j}\}.
	\end{equation}
	We define the representation \eqref{5.7x} by taking
	\begin{equation}\label{10.5}
\chi(\gamma_i)=E_i^{j}.
	\end{equation}

To verify \textbf{[BIV]} \eqref{5.10}, let $\gamma=\alpha_m\cdots \alpha_1$ be a nonempty reduced word, where each $\alpha$ is one of $\gamma_i^\pm$, then by \eqref{10.4x},
\begin{equation}\label{10.6}
\chi(\gamma)\cdot\mathbb{R}\begin{pmatrix}0\\1\end{pmatrix}\in D_{\chi(\alpha_m)},\qquad \mathbb{R}\begin{pmatrix}0\\1\end{pmatrix}\notin D_{\chi(\alpha_m)}.
\end{equation}
In particular, we have $\chi(\gamma)\neq I$, proving that $\chi(\Gamma)$ is freely generated. For a nonidentity $\chi(\gamma)=\blk\begin{smallmatrix}a_\gamma&b_\gamma\\c_\gamma&d_\gamma\end{smallmatrix}\brk$, the condition $b_\gamma=0$ would imply $\chi(\gamma)\cdot \mathbb{R}\big(\begin{smallmatrix}0\\1\end{smallmatrix}\big)=\mathbb{R}\big(\begin{smallmatrix}0\\1\end{smallmatrix}\big)$, contradicting \eqref{10.6}.

To verify \textbf{[ERG]} \eqref{5.11}, note that every generator $\chi(\gamma_i)=E_i^j$ is hyperbolic and hence acts ergodically on $\mathbb{T}^2$. This is stronger than \textbf{[ERG]}, which requires only that functions invariant under the whole group $\chi(\Gamma)$ be constant almost everywhere. Here, invariance under any single generator already suffices.

Having shown that both \textbf{[BIV]} \eqref{5.10} and \textbf{[ERG]} \eqref{5.11} can hold simultaneously, we now briefly discuss how they can fail. \textbf{[BIV]} fails if some nonidentity $\gamma\in\Gamma$ has a singular $B_\gamma$-block. For example, if $\chi(\gamma)=I$ for some nonidentity $\gamma\in\Gamma$. \textbf{[ERG]} fails if some nonzero $\omega\in\mathbb{Z}^{2n}$ has a finite orbit under $\chi(\Gamma)^T$. For example, if every generator is mapped to an upper triangular matrix $\blk\begin{smallmatrix}1&a_i\\0&1\end{smallmatrix}\brk$, then $\big(\begin{smallmatrix}0\\1\end{smallmatrix}\big)$ is fixed by every transpose. Hence the second coordinate on $\mathbb T^2$ is invariant, and the action is not ergodic.

	\addcontentsline{toc}{section}{References}
	\bibliographystyle{alpha}
	\def\cprime{$'$} \def\cprime{$'$}

\end{document}